\documentclass[a4paper,12pt]{amsart}

\usepackage{amsmath,amssymb}
\usepackage[T1]{fontenc}
\usepackage{latexsym}

\def \ni{\noindent}
\def \vs{\vskip .6cm}

\def \t{\widetilde}
\def \beq{\begin{eqnarray*}}
\def \eeq{\end{eqnarray*}}

\def \n{\nabla}
\def \nt{\widetilde{\nabla}}
\def \nb{\overline \nabla}
\def \Rt{\widetilde R}
\def \Rb{\overline R}

\def \h{\mathfrak h}
\def \g{\mathfrak g}
\def \k{\mathfrak k}
\def \v{\mathfrak v}
\def \l{\mathfrak l}
\def \m{\mathfrak m}

\def \s{\mathfrak s}
\def \t{\mathfrak t}
\def \n{\mathfrak n}
\def \A{\mathcal A}
\def \a{\mathfrak a}
\def \H{\mathcal H}
\def \V{\mathcal V}
\def \p{\mathfrak p}

\def \LLambda {{\bf \Lambda}}
\def \llambda {\boldsymbol{\lambda}}
\def \so {\boldsymbol{\mathfrak{so}}}
\def \u {\boldsymbol{\mathfrak u}}
\def \Z {{\bf Z}}
\def \ZZ {\boldsymbol{\mathcal Z}}

\def \RM{\mathbb{R}}
\def \NM{\mathbb{N}}
\def \ZM{\mathbb{Z}}
\def \CM{\mathbb{C}}
\def \HM{\mathbb{H}}

\def \J {\mathcal J}
\def \Jt {\widetilde{\mathcal J}}
\def \Jb {\overline{\mathcal J}}
\def \N {\mathcal N}
\def \Nt {\widetilde{\mathcal N}}

\newcommand \bo [1] {\bigotimes \phantom{ }^#1}

\def \leftr {[\hbox{\hspace{-0.15em}}[}
\def \rightr {]\hbox{\hspace{-0.15em}}]}
\def \la {\langle}
\def \ra {\rangle}
\def \bo {\bigoplus}

\newtheorem{defi}{Definition}[section]

\newtheorem{prop}[defi]{Proposition}

\newtheorem{theo}[defi]{Theorem}

\newtheorem{lemm}[defi]{Lemma}

\newtheorem{NB}[defi]{Remark}

\newtheorem{coro}[defi]{Corollary}

\begin{document}

\begin{center}
{\Large \bf Twistors and 3-symmetric spaces}
\vs

{\large Jean-Baptiste Butruille} \vs

\end{center}

\vs

{\small 

\ni{\bf Abstract} --  We describe complex twistor spaces over inner 3-symmetric spaces $G/H$, such that $H$ acts transitively on the fibre. Like in the symmetric case, these are flag manifolds $G/K$ where $K$ is the centralizer of a torus in $G$. Moreover, they carry an almost complex structure defined using the horizontal distribution of the normal connection on $G/H$, that coincides with the complex structure associated to a parabolic subgroup $P \subset G^{\mathbb C}$ if it is integrable. Conversely, starting from a complex flag manifold $G^{\mathbb C}/P$, there exists a natural fibration with complex fibres on a 3-symmetric space, called fibration of degree 3.


}

\vs \vs

\section{Introduction}

Twistor theory is a story that begins in dimension 4 : let $(M,g)$ be an oriented Riemannian 4-manifold, because of the decomposition, at the Lie algebra level, $\mathfrak{so}(4) \simeq \mathfrak{sp}(1) \oplus \mathfrak{sp}(1)$, the bundle of 2-forms splits into
\[ \LLambda^2 = \LLambda^2_+ \oplus \LLambda^2_- \]
The twistor space of $M$ is defined to be $\ZZ$, the unit sphere bundle in $\LLambda^2_-$. It can be seen as the bundle whose sections are the almost complex structures (or the Kähler forms) on $M$, compatible with the metric and orientation. It is this last definition we use in even dimension higher than 6. The fibre $\ZZ_x \simeq \CM P(1)$, or the vertical distribution, are equipped with a natural, and in fact with two opposite complex structures which may be completed using the horizontal distribution induced by the Levi Civita connection : the resulting almost complex structure $J_1$ on $\ZZ$ is integrable, for one choice of orientation of the fibre, as soon as the manifold is self-dual, meaning that the Weyl tensor takes values in $\LLambda^2_+ $, while with the other choice of orientation, we get an almost complex structure $J_2$ never integrable.

Unfortunately, this has much to do with the representations of $SO(4)$ and generalizes quite bad in higher dimensions. However, solutions to the problem of finding complex twistor spaces over a Riemannian manifold $M$ of dimension $2n \geq 6$, non locally conformally flat, exist if we suppose it is equipped furthermore with some {\it $G$-structure}, $G \varsubsetneq SO(2n)$. We enter then the field of {\it special geometries}. So far, the research on special geometries and twistor spaces has been mainly focused, quite naturally, on special holonomies (for example Kähler, or quaternion-Kähler manifolds), including symmetric spaces. We are interested here in special geometries {\it with torsion}, i.e. the case where the $G$-structure is not underlying a Riemannian holonomy reduction.

Nearly Kähler manifolds are quite representative of the special geometries with torsion, for several reasons. Firstly they consist in almost Hermitian manifolds $(M,g,J)$ which are not Kählerian, but not either complex ($J$ is not integrable), nor symplectic (the Kähler form $\omega$ is not closed). So they are as far as possible of integrable structures. However, the condition imposed on the $G$-structure to define them is quite simple : we ask $\nabla \omega$, where $\nabla$ is the Levi-Civita connection, to be a 3-form, $\forall X \in TM$, $\nabla_X \omega = \frac{1}{3}\iota_X d\omega$. Secondly, the classification of nearly Kähler manifolds is intimately related to a class of homogeneous, non symmetric manifolds (such a homogeneous space $K/H$ is thus an example of a {\it $H$-manifold with torsion}) : the 3-symmetric spaces. Indeed, it was shown in \cite{bu} that all homogeneous nearly Kähler manifolds are of that type and in fact these are the only compact (or equivalently complete) known examples. And thirdly, in dimension 6, nearly Kähler manifolds have a further reduction to $SU(3)$ and belong, though not to the series of special holonomies, to another series : the weak holonomies, like Einstein-Sasaki manifolds in dimension 5 or nearly parallel $G_2$-manifolds in dimension 7. Kähler geometry, Riemannian holonomy, symmetric spaces : the study of nearly Kähler manifolds gives a generalization of these three notions that doesn't consist only in weakening the definition but in imposing new, orthogonal conditions that still determine useful identities, for example on the curvature (this constitutes one of the motivations given by the major contributor of the field, Alfred Gray, see e.g. \cite{gr2}).

Moreover, nearly Kähler manifolds are already related to twistor theory. Indeed, it was shown by Eels, Salamon \cite{ee} that the twistor space $\ZZ$ of a self-dual Einstein manifold has a natural Kähler structure $(g_1,J_1)$ (cf \cite{fr,hi2}) but also a natural nearly Kähler structure $(g_2,J_2)$. The same holds, as proved in \cite{al} or \cite{na}, on the twistor space of a quaternion-Kähler manifold with positive scalar curvature. And if the base is symmetric, then $\ZZ$ is a 3-symmetric space. On the other hand, $S^6 \simeq G_2/SU(3)$ has a complex Riemannian twistor space, because it is locally conformally flat, and inside this one, a {\it reduced twistor space}, isomorphic to $G_2/U(2)$, invariant by $SU(3)$.

Inspired by this last example, we look, in this article like in a previous one \cite{bu2}, for twistor spaces with integrable almost complex structure on nearly Kähler manifolds. According to Nagy \cite{na2}, nearly Kähler, complete, simply connected manifolds decompose, in all dimensions, in a Riemannian product of : \\
(\romannumeral 1) 3-symmetric spaces, \\
(\romannumeral 2) twistor spaces of non locally symmetric, irreducible, quaternion-Kähler manifolds, \\
(\romannumeral 3) 6-dimensional, non locally 3-symmetric, strictly nearly Kähler manifolds.

In \cite{bu2} we were concerned by (\romannumeral 3), i.e. 6-dimensional manifolds, going back to a problem left opened by O'Brian, Rawnsley \cite{ob}. The results were not fully satisfactory, from our point of view, since only $S^6$ and Kählerian manifolds, or almost Hermitian manifolds locally conformal to them, were given a {\it complex} twistor space with our method. We consider now the case of (\romannumeral 1). This time we were able to generalize the fibration $G_2/U(2) \to G_2/SU(3) \simeq S^6$. The reason is we only used in \cite{bu2} the reduction of the structural group to $U(n)$ while we take advantage here of the full holonomy reduction $H$ of $\nb$, the canonical Hermitian connection of our 3-symmetric space $G/H$.

The paper is organized as follows. Section 2 is devoted to some quite general preliminaries. We give the definition of nearly Kähler manifolds. Section 3 is a presentation of 3-symmetric spaces based on the isotropy representation. We show that the normal connection of a Riemannian 3-symmetric space (3-symmetric homogeneous spaces are always {\it reductive}) coïncides with the intrinsic connection of a canonical almost Hermitian structure defined on it. As a consequence, we are able to give a new characterization of locally 3-symmetric spaces (see theorem \ref{locally 3-symmetric}) as a particular case of Ambrose-Singer manifolds. However, it can be seen as a reinterpration, in terms of the intrinsic connection, of equations given by Gray in \cite{gr2}. Moreover, it generalizes a remark of Nagy in \cite{na2} in the nearly Kähler case. In section 4 we explain how to construct a twistor space $\Z$ and its almost complex structure on a $H$-manifold with torsion, such that the group $H$ acts transitively on the fibre, using a $H$-connection. In section 5 we extend the work of Burstall, Rawnsley \cite{bur} and others on symmetric spaces to a general reductive homogeneous space $G/H$, such that the groups $G$ and $H$ have same rank. We are particularly interested in the case of 3-symmetric spaces (then, our hypothesis is equivalent to $s$, the automorphism of order 3 defining our space, being {\it inner}, like in the symmetric case). Twistor spaces constructed on $G/H$, following the plan of section 4, are flag manifolds. By this we mean a homogeneous space $G/K$ where $K$ is the centralizer of a torus in $G$. In particular $G$ acts transitively on them and integrability conditions (\ref{Jt integrable torsion}), (\ref{Jt integrable curvature}), translate into algebraic conditions (\ref{Jt integrable homogene}) on a subspace $\n^+$ of the complexified Lie algebra $\g^{\CM}$. Invariant complex structures on flag manifolds are given by parabolic subgroups $P \subset G^{\CM}$ such that $K = P \cap G$ and so $G/K \simeq G^{\CM}/P$. Our main result of this section (theorem \ref{G^c/P -> G/H}) is then that the natural almost complex structure associated to the normal connection on the twistor space is integrable if and only it is given by such an isomorphism. We summarize this :
\begin{theo}
Let $M=G/H$ be a normal homogeneous space such that $G$ and $H$ have same rank. Let $(\ZZ,\Jb)$ be the Riemannian twistor space of $M$, equipped with the almost complex structure associated to the normal connection. Then \\
(\romannumeral 1) $\Z$ is an \emph{almost complex} submanifold of $(\ZZ,\Jb)$ if and only if it is isomorphic to a flag manifold $G/K$, where $K$ is a subgroup of $H$. \\
(\romannumeral 2) Supposing (\romannumeral 1), $(\Z,\Jb)$ is a \emph{complex} manifold if and only if it is isomorphic to a complex flag manifold $G^{\CM}/P$.
\end{theo}
Conversely, starting with a complex flag manifold $G^{\CM}/P$, we define a fibration with complex fibres over a 3-symmetric space, called fibration of degree 3. We use the same tools as Burstall and Rawnsley \cite{bur} for the construction of their {\it canonical fibration}, going from a flag manifold to a {\it symmetric} space. The results of section 5 are applied, in the last two sections, to 3-symmetric spaces : in section 6, to examples which we found enlightening, and in section 7 systematically to two types of inner 3-symmetric spaces in the classification of Gray and Wolf \cite{wo} : isotropy irreducible 3-symmetric spaces, like the sphere $S^6 \simeq G_2/SU(3)$ ; and twistor spaces over symmetric spaces, like $\CM P(2q+1)$, the twistor space of $\HM(q)$. The former were first discovered by Wolf in \cite{wo3} : indeed, in the list given p281, all homogeneous spaces are symmetric, except six 3-symmetric spaces and one 5-symmetric space.

\begin{theo}
Isotropy irreducible 3-symmetric spaces $G/H$ have a complex twistor space $G/K$ such that \\
(\romannumeral 1) the fibre, $H/K$, is isomorphic to $\CM P(2) \simeq SU(3)/U(2)$, except for $E_8/SU(9)$ where the fibre is $\CM P(8) \simeq SU(9)/U(8)$. \\
(\romannumeral 2) $G/K \to G/H$ is the fibration of degree 3.
\end{theo}

\ni {\it Acknowledgement :} thank you as always to Andrei (Moroianu) for his precious advice.

\section{Preliminaries}

The objects -- nearly Kähler, 3-symmetric spaces -- we deal with in this article are both "algebraic" (since they are isomorphic to a quotient of Lie groups) and "geometric" (since they are almost Hermitian differentiable manifolds). So we need two sorts of preliminaries. First, concerning Lie algebras and representation theory :

\begin{defi}
Let $\a$ be a Lie algebra. The derived series $(\a^{(i)})_{i \geq 1}$ of $\a$ is defined by
\[ 
\a^{(1)} = \a \quad \text{and} \quad \a^{(i+1)} = [\a^{(i)},\a^{(i)}]
\]
The central descending series of $\a$ are $(\a^i)_{i \geq 1}$ where
\[ 
\a^1 = \a, \qquad \a^{i+1} = [\a^1,\a^i]
\]
\end{defi}

\begin{defi}
A Lie algebra $\a$ is called solvable if there exists an $r \in \NM$ such that $\a^{(r)}=\{0\}$. It is called nilpotent if there exists $r$ such that $\a^r=\{0\}$.
\end{defi}

\begin{defi}
Let $\g$ be the Lie algebra of a compact semi-simple Lie group $G$. \\
A parabolic subalgebra $\p$ is a subalgebra that contains a maximal solvable subalgebra of $\g^{\CM}$.
\end{defi} 

The parabolic subalgebras are described using a root system. Let $T \subset G$ be a maximal torus, $\t$ its Lie algebra. We denote by $\mathcal R$ the root system associated to $T$, i.e. the set of non zero weights of the adjoint representation ${\rm Ad}$ of $G$ on the complexified Lie algebra $\g^{\CM} = \g \otimes_{\RM} \CM$. Let $\mathcal B$ be a choice of simple roots. For $\beta \in \mathcal B$, $\alpha \in \mathcal R$, denote by $n_{\beta}(\alpha)$ the coefficient along $\beta$ of $\alpha$ in the decomposition
\[ \alpha = \sum_{\beta \in \mathcal B} n_{\beta}(\alpha) \beta \]
relative to the base $\mathcal B$. By the definition of the simple roots, $\alpha$ is a positive root ($\alpha \in \mathcal R^+$) if and only if $n_{\beta}(\alpha) > 0$ {\it for all} $\beta \in \mathcal B$. More generally define, for a subset $\A \subset \mathcal B$,
\[ n_{\A} = \sum_{\beta \in \mathcal A} n_{\beta} \]
Again, $\alpha \in \mathcal R^+$ implies $n_{\A}(\alpha) \geq 0$ and $n_{\A}(\alpha)=0$ if and only if $n_{\beta}(\alpha)=0$ for all $\beta \in \A$. 

Denote by $\g_{\alpha} \subset \g^{\CM}$ the 2-dimensional weight space relative to a root $\alpha$. The weight space relative to 0 is simply $\t^{\CM}$. Denote, also, by $\leftr \g_{\alpha} \rightr$ the {\it real} root space :
\[ \leftr \g_{\alpha} \rightr = \g \cap (\g_{\alpha} \oplus \g_{-\alpha}). \]

\begin{theo}
For each subset $\A$ of $\mathcal B$, the summand
\begin{equation}
\p_{\A} = \t^{\CM} \oplus \sum_{n_{\A}(\alpha) \geq 0} \g_{\alpha}
\label{pA}
\end{equation}
is a parabolic subalgebra of $\g^{\CM}$ and moreover, each parabolic subalgebra can be written in this form for a certain root system and a choice of simple roots.
\label{parabolic subalgebras}
\end{theo}

We decompose $\p_{\A}$ into
\begin{equation}
\p_{\A} = \k_{\A}^{\CM} \oplus \l_{\A}^+ 
\label{p=k^c + l-}
\end{equation}
where
\begin{equation}
\k_{\A} = \g \cap \p_{\A} = \t \oplus \sum_{n_{\mathcal A}(\alpha) = 0} \leftr \g^{\alpha} \rightr, \qquad \l_{\A}^+ = \sum_{n_{\A}(\alpha) > 0} \g^{\alpha} 
\label{kA,lA+}
\end{equation}
Let
\[ \l_{\A}^- = \overline{\l_{\A}^+} = \sum_{n_{\A}(\alpha) < 0} \g^{\alpha}, \]
\[ \g^{\CM} = \p_{\A} \oplus \l_{\A}^- \]

\begin{lemm}
The summand $\l_{\A}^+$ is a nilpotent subalgebra. In fact, it is the nilradical (the maximum nilpotent subalgebra) of $\p_{\A}$.
\label{lA+ nilradical}
\end{lemm}

\vs

Let's come to the geometric preliminaries. They deal with almost Hermitian manifolds, then with nearly Kähler manifolds.

An almost Hermitian manifold is a real manifold $M$ of dimension $m=2n$ with a reduction of the frame bundle to the unitary group $U(n)$ (a {\it $U(n)$-structure}) or equivalently with a Riemannian metric $g$ and an almost complex structure $J$ ($J^2=-Id$), compatible in the sense that, $\forall X,Y \in TM$, $g(JX,JY)=g(X,Y)$. Given $g$, $J$, we construct a third tensor, a differential 2-form $\omega$, called the Kähler form of $M$ and defined by
\[ \omega(X,Y) = g(JX,Y). \]
Let $\nabla$ be the Levi-Civita connection of $g$. The set of metric connections is an affine space modeled on the space of sections of $\LLambda^1 \otimes \so(m)$, where $\so(m)$ is the bundle of skew-symmetric endomorphisms of $TM$. In other words, if $\nt$ is a metric connection, the difference $\nabla - \nt$ is a $\so(m)$-valued 1-form. Now the set of Hermitian connections, that is connections $\nt$ satisfying not only $\nt g = 0$ but also $\nt J = 0$, is an affine subspace of the previous one, modeled on $\Gamma(\LLambda^1 \otimes \u(n))$, where $\u(n)$ is the adjoint bundle of the $U(n)$-structure, identified with the bundle of skew-symmetric endomorphisms of the tangent bundle commuting with $J$. In other words, $\nabla - \nt$ decomposes into
\begin{equation}
\begin{array}{ccccc}
\LLambda^1 \otimes \so(m) & = & \LLambda^1 \otimes
\u(n)^{\perp} & \oplus & \LLambda^1 \otimes \u(n) \\
\nabla - \nt & = & \eta & + & \xi
\end{array}
\label{n-nt}
\end{equation}
where the tensor $\eta$ is independant on the choice of the Hermitian connection. Let $\nb$ be the unique Hermitian connection such that $\xi = 0$ or
\[ \nabla - \nb = \eta \]

\begin{defi}
We call $\nb$ the intrinsic connection, or the canonical Hermitian connection, of $(M,g,J)$. Then we call $\eta$ the intrinsic torsion of the $U(n)$-structure or the almost Hermitian manifold.
\end{defi}

The torsion $T$ of $\nb$ and $\eta$, the intrinsic torsion, are related by
\[ T(X,Y) = \eta_X Y - \eta_Y X \]
Furthermore $\eta$ can be computed explicitely :
\[ \forall X \in TM, \quad \eta_X = \frac{1}{2} J \circ (\nabla_X J) \]
Then, $\eta=0$ if and only if $J$ is parallel for the Levi-Civita connection, i.e. $M$ is Kählerian, or the Riemannian holonomy is contained in $U(n)$. Now,

\begin{defi}
An almost Hermitian manifold $(M,g,J)$ is called nearly Kähler if the intrinsic torsion is totally skewsymmetric,
\[ \forall X \in TM, \quad (\nabla_X J) X = 0 \]
Equivalently, the covariant derivative of $\omega$ is a 3-form :
\[ \nabla \omega = \frac{1}{3} d\omega \]
because the Levi-Civita connection has no torsion.
\end{defi}

Denote by $\llambda^{p,q}$ the bundle of $r$-forms of type $(p,q)$, $p+q=r$. If $p \neq q$, the intersection of $\llambda^{p,q} \oplus \llambda^{q,p}$ with the real forms is the bundle $\leftr \llambda^{p,q} \rightr \subset \LLambda^r$ of real forms of type $(p,q)+(q,p)$. But if $p=q$, $\llambda^{p,p}$ is already the complexification of $[\llambda^{p,p}]$, the bundle of real forms of type $(p,p)$. Now, the metric $g$ gives the isomorphisms $\u(n) \simeq [\llambda^{1,1}]$ and $\u(n)^{\perp} \simeq \leftr \llambda^{2,0} \rightr$ and there is a decomposition,
\begin{eqnarray} 
\LLambda^1 \otimes \u(n)^{\perp} & \simeq & \leftr \llambda^{1,0} \otimes \llambda^{2,0} \rightr \oplus \leftr \llambda^{0,1} \otimes \llambda^{2,0} \rightr \label{W12+34} \\
& \simeq & \leftr \llambda^{3,0} \rightr \oplus \leftr {\bf U^1} \rightr \oplus \leftr
\llambda^{2,1}_0 \rightr \oplus \LLambda^1
\end{eqnarray}
This is the irreducible decomposition of the representation $\LLambda^1 \otimes \u(n)^{\perp}$ of $U(n)$ at each point. The totally skewsymmetric tensors are the sections of the bundle isomorphic to $\leftr \llambda^{3,0} \rightr$. Then, $M$ is nearly Kähler if and only if $\eta \in \leftr \llambda^{3,0} \rightr$ whereas complex manifolds are characterized by $\eta \in \leftr
\llambda^{2,1} \rightr$, and $d\omega = 0$ (the manifold $M$ is symplectic, or {\it almost Kähler}) is equivalent to $\eta$ being a section of the unidentified bundle $\leftr {\bf U^1} \rightr$.

With this definition, Kählerian manifolds are nearly Kähler but we want to treat them separately so we give the following definition :
\begin{defi}
A nearly Kähler manifold is called strict (strictly nearly Kähler manifold) if there exist no vector $X \in TM$ such that $\nabla_X J$ is identically zero.
\end{defi}

Nearly Kähler manifold have the following important property :
\begin{equation} 
\nb \eta = 0 
\label{intrinsic torsion //}
\end{equation}
Since $\nb$ is a Hermitian connection, this is equivalent to $\nb (\nabla \omega) = 0$ or $\nb d\omega = 0$.

A consequence (cf \cite{na}) is that any nearly Kähler manifold might be decomposed locally (or globally, if it is  simply connected) into the Hermitian product of a Kähler manifold by a {\it strictly} nearly Kähler manifold. Thus, the study of nearly Kähler manifolds reduces to that of strict ones. From now on, when talking about nearly Kähler manifolds we implicitely assume they are strict.

We end up this section by quoting a result announced in the introduction. This is due to Alexandrov, Grantcharov, Ivanov \cite{al} or Nagy \cite{na} for quaternion-Kähler manifolds, and to Eels, Salamon \cite{ee} in dimension 4.
\begin{theo}
Let $\pi : M \to N$ be the twistor space of a quaternion-Kähler manifold of dimension $4q$, $q \geq 2$, or of a Kähler-Einstein manifold in dimension 4, with positive scalar curvature. Then $M$ has a natural Kählerian structure  $(g_1,J_1)$ (see \cite{hi2,fr,sa3}) but also a natural strictly nearly Kähler structure $(g_2,J_2)$. Moreover, let $\V$ be the vertical distribution on $M$, tangent to the fibres, and let $\H$ be the horizontal distribution (such that $TM = \V \oplus \H$) induced by the Levi-Civita connection of $g$, the metric on $N$, we have :
\[ g_1|_{\V} = 2 g_2|_{\V}, \quad g_1|_{\H} = g_2|_{\H} = \pi^*(g) \quad \text{and} \quad J_1|_{\V} = - J_2|_{\V}, \quad J_1|_{\H} = J_2|_{\H} \]
\label{QK twistor is NK}
\end{theo}

\section{3-symmetric spaces}

Another important class of examples of nearly Kähler manifolds is provided by the 3-symmetric spaces. They were introduced by Gray \cite{gr2} in 1970 as a generalization of the well-known symmetric spaces :

\begin{defi}
A 3-symmetric space is a homogeneous space $M=G/H$, where $G$ has an automorphism $s$ of order 3 such that
\begin{equation} G^s_0 \subset H \subset G^s \label{H fixed points} \end{equation}
where $G^s=\{ g \in G \ | \ s(g)=g \}$ is the fixed points set of $s$ and $G^s_0$ is the identity component of $G^s$.
\end{defi}

Let $\g$, $\h$ be the Lie algebras of $G$, $H$, respectively. For a symmetric space, the eigenspace, for the eigenvalue $-1$, of the differential $s_* : \g \to \g$ of the automorphism at $e$, is an ${\rm Ad}(H)$-invariant complement of $\h$ in $\g$, so that a symmetric space is always a {\it reductive} homogeneous space. 

For a 3-symmetric space, things are slightly more complicated : $s_*$ has now three eigenvalues, $1$, $j=-\frac{1}{2} + i\frac{\sqrt{3}}{2}$ and $j^2=\bar j = -\frac{1}{2} - i\frac{\sqrt{3}}{2}$, two of which are complex, and correspondingly there is a decomposition (by (\ref{H fixed points}), $\h$ is always the eigenspace for the eigenvalue 1) :
\[ g^{\CM} = \h^{\CM} \oplus \m_j \oplus \m_{j^2} \]
Now, $\m = (\m_j \oplus \m_{j^2}) \cap \g$ satisfies 
\begin{equation} \g = \h \oplus \m, \quad {\rm Ad}(H)\m \subset \m \label{reductive} \end{equation}
so 3-symmetric spaces are also reductive homogenous spaces. 

But we have more than this. 

Recall that on a reductive homogenous space, the associated bundle $G \times_{{\rm Ad}} \m$ (we view $G$ as a principal bundle of structure group $H$ over $M$), is identified with the tangent bundle $TM$ and that invariant tensors, for the right action of $G$ on $M$, are identified with constant tensors on $\m$. Recall also that an almost complex structure might be identified with the bundle $T^+M$ of (1,0) vectors. Thus, an {\it invariant} almost complex structure on a reductive homogeneous space is identified with a maximal isotropic subspace $\m^+$ of $\m^{\CM}$ or equivalently with a decomposition $\m^{\CM} = \m^+ \oplus \m^-$ where $\m^- = \overline{\m^+}$. Now, on a 3-symmetric space there is already such a decomposition, by the definition of $\m$.
\begin{defi}
The canonical almost complex structure of a 3-symmetric space is the invariant almost complex structure such that the (1,0)-vector fields are the sections of the associated bundle $G \times_{{\rm Ad}} \m_j$. 
\end{defi}
In other words we put $\m^+ = \m_j$. We will give another definition, more faithful to the original one. The restriction of $s_*$ to $\m$ represents an invariant tensor $S$ of $M$ which satisfies : \\
$(\romannumeral 1)$ $S^3 = Id$ \\
$(\romannumeral 2)$ $\forall x \in M,$ 1 is not an eigenvalue of $S_x$ \\
Thus, one can write $S$, as for a (non trivial) third root of unity,
\begin{equation} 
S = -\frac{1}{2}Id + \frac{\sqrt{3}}{2}J 
\label{canonical J}
\end{equation}
where
\[ J^2 = -Id \]
We say that $J$ is the canonical almost complex structure of the 3-symmetric space.

Similarly, on a reductive homogeneous space, an ${\rm Ad}(H)$-invariant scalar product $g$ on $\m$, where $\m$ satisfies (\ref{reductive}), defines an invariant metric on $M$, also denoted by $g$. If $M$ is a symmetric space, the pair $(M,g)$ is called a {\it Riemannian} symmetric space. If $M$ is a 3-symmetric space and $g$, $J$ are compatible, $(M,g)$ is called a Riemannian 3-symmetric space.

A sufficient condition for a reductive homogeneous space to be a symmetric space is
\begin{equation} 
[\m,\m] \subset \h \label{symmetric} 
\end{equation}
Indeed, we define an endormorphism $f$ of $\g$ by setting $f|_{\h} = Id|_{\h}$ and $f|_{\m} = -Id|_{\m}$. Then, (\ref{symmetric}) together with
\begin{equation} 
[\h,\m] \subset \m \label{reductive inf} 
\end{equation}
(a consequence of (\ref{reductive})) imply
\begin{equation} 
\forall X,Y \in \g, \quad [f(X),f(Y)]=f([X,Y]) 
\label{g automorphism} 
\end{equation}
So $f$ can be integrated into an automorphism (an involution) of $G$ :
\[ \forall X \in \g, \quad s(\exp(X))=\exp(f(X)) \]
($G$ is simply connected whence $\exp : \g \to G$ is surjective).

Similar conditions exist, for 3-symmetric spaces :
\begin{lemm}
Let $M = G/H$ be a reductive homogeneous space, with an invariant almost complex structure represented by the decomposition
\[ \m^{\CM} = \m^+ \oplus \m^- \]
Then $M$ is a 3-symmetric space endowed with its canonical almost complex structure if and only if
\begin{equation} 
[ \m^+,\m^+ ] \subset \m^-, \quad [ \m^-,\m^- ] \subset \m^+ \quad \text{and} \quad [ \m^+,\m^- ] \subset \h^{\CM}.
\label{3symmetric} 
\end{equation}
\end{lemm}

\begin{proof}
The proof is very similar to the symmetric case. We set
\[ f|_{\h} = Id|_{\h}, \quad f|_{\m^+} = j Id|_{\m^+} \quad \text{and} \quad f|_{\m^-} = j^2 Id|_{\m^-} \]
The endomorphism $f$ satisfies (\ref{g automorphism}) if and only if, for two eigenspaces $\m_{\lambda}$ and $\m_{\mu}$ corresponding to the eigenvalues $\lambda$, $\mu$ we have $[\m_{\lambda},\m_{\mu}] \subset \m_{\lambda \mu}$. We already have : $[\h,\h] \subset \h$ ($\h$ is a subalgebra) and
\[
[ \h,\m^+ ] \subset \m^+ \quad \text{and} \quad [ \h,\m^+ ] \subset \m^+,
\]
because the almost complex structure represented by $\m^+$ is invariant. The rest is exactly (\ref{3symmetric}).
\end{proof}

Conditions involving the Lie bracket might be interpreted, on a reductive homogeneous space, as conditions on the torsion and the curvature of the normal connection $\widehat \nabla$. The latter is defined as the $H$-connection on $G$ whose horizontal distribution is $G \times \m \subset TG \simeq G \times \g$. Indeed,

\begin{lemm}
The torsion $\widehat T$, and the curvature $\widehat R$ of the normal connection $\widehat \nabla$, viewed as constant tensors, are respectively the $\m$-valued 2-form and the $\h$-valued 2-form on $\m$ given by
\[ \forall u,v \in \m, \quad \widehat T(u,v)=-[u,v]^{\m}, \quad \widehat R_{u,v} = [u,v]^\h \]
\label{normal torsion curvature}
\end{lemm}

\begin{prop}
Let $M = G/H$ be an almost Hermitian homogeneous space. By this we mean it is equipped with a $G$-invariant almost Hermitian structure. Suppose it is furthermore reductive. Then, $M$ is a 3-symmetric space if and only if it is quasi-Kähler and the intrinsic connection $\nb$ coïncides with $\widehat \nabla$.
\label{connexion normale = intrinseque}
\end{prop}

An almost Hermitian manifold is called {\it quasi-Kähler} if $\forall X,Y \in TM$, $(\nabla_X J)Y + (\nabla_{JX} J)JY=0$. Equivalently, the intrinsic torsion is a section of the first bundle of (\ref{W12+34}), isomorphic to $\leftr \llambda^{1,0} \otimes \llambda^{2,0} \rightr$. 

\begin{proof}
In terms of the torsion and curvature of the normal connection, equations (\ref{3symmetric}) become
\begin{equation} 
\widehat T(\m^+,\m^+) \subset \m^-, \quad \widehat T(\m^-,\m^-) \subset \m^+, \quad \widehat T(\m^+,\m^-) = \{0\}
\label{3symmetric torsion}
\end{equation}
\begin{equation}
\widehat R(\m^+,\m^+) = \widehat R(\m^-,\m^-) = \{ 0 \}
\label{3symmetric courbure}
\end{equation}

The first line tells us that $\widehat T$ belongs to the intersection of $\Lambda^2 \m \otimes \m$ with $\bigotimes^3 (\m^+)^* \oplus \bigotimes^3 (\m^-)^*$. Since the application
\begin{eqnarray*}
\LLambda^1 \otimes \so(m) & \to & \LLambda^2 \otimes TM \\
\widetilde \eta & \mapsto & \{ \widetilde T : (X,Y) \mapsto \widetilde T(X,Y)=\widetilde \eta_X Y - \widetilde \eta_Y X \}
\end{eqnarray*}
is an isomorphism, we prefer working with $\widehat \eta = \nabla - \widehat \nabla$. Then (\ref{3symmetric torsion}) is equivalent to
\[ \widehat \eta \in \leftr \llambda^{1,0} \otimes \llambda^{2,0} \rightr \]
In particular $\widehat \eta \in \LLambda^1 \otimes \u(n)^{\perp}$ so $\widehat \nabla = \nb$ (the normal connection of an almost Hermitian homogeneous space is always a Hermitian connection) thus $\eta \in \leftr \llambda^{1,0} \otimes \llambda^{2,0} \rightr$, i.e. the manifold is quasi-Kähler and the converse is true (equations (\ref{3symmetric courbure}) are automatically satisfied on a quasi-Kähler manifold by the curvature tensor of the intrinsic connection, see for example \cite{fa}).
\end{proof}

There is also a notion of {\it locally 3-symmetric space}, which means that there exists a family of local isometries $(s_x)_{x \in M}$, the {\it geodesic symmetries of order 3}, such that $s^3_x = Id$ and $\forall x \in M$, $x$ is an isolated fixed point of $s_x$ (for a 3-symmetric space, the automorphism $s$ induces such a family on the manifold where the $s_x$ are globally defined). See \cite{gr2} for details.

\begin{theo}
An almost Hermitian manifold $M$ of dimension $m$ is a locally 3-symmetric space if and only if it is quasi-Kähler and the torsion and the curvature of the intrinsic connection $\nb$ satisfy
\begin{equation} 
\nb T = 0 \quad \text{and} \quad \nb \, \Rb = 0 
\label{T, Rb paralleles}
\end{equation}
\label{locally 3-symmetric}
\end{theo}

\begin{proof}
Let $H$ be the reduced holonomy group of the intrinsic torsion and $\h$ its Lie algebra. If the torsion and curvature of $\nb$ are parallel, they generate an infinitesimal model $(\g,\h)$. Moreover, it is always {\it regular}, like in the symmetric case, under the hypothesis of theorem (\ref{locally 3-symmetric}) because $\h$ is the fixed point set of a Lie algebra automorphism of order 3 : if $G$ be the simply connected group of Lie algebra $\g$, $H$ is a closed subgroup of $G$. Then, by the work of Tricerri \cite{tr2}, $M$ is locally isometric to a reductive homogeneous space $G/H$ whose normal connection coïncides with $\nb$ and since $M$ is quasi-Kähler, by proposition \ref{connexion normale = intrinseque}, $G/H$ is a 3-symmetric space.

Conversely, the normal connection of a reductive homogeneous space always satisfies $\widehat \nabla \widehat T = \widehat \nabla \widehat R = 0$ (in fact, every invariant tensor for the action of $G$ on the right is parallel for $\widehat \nabla$) and on a 3-symmetric space it coïncides with $\nb$.
\end{proof}

The system (\ref{T, Rb paralleles}) should be compared to the definition, $\nabla R = 0$, of locally symmetric spaces. In \cite{bu3} we showed that it is equivalent to some equations given by Gray \cite{gr2}, reinterpreted by the intrinsic connection.

Nearly Kähler manifolds are quasi-Kähler and the condition that $M$, the 3-symmetric space, is nearly Kähler translates into a structural condition on the homogeneous space.

\begin{defi}
A reductive Riemannian homogeneous space is called naturally reductive if the scalar product $g$ on $\m$ representing the invariant metric satisfies
\begin{equation} 
\forall X,Y,Z \in \m, \quad g([X,Y],Z)=-g([X,Z],Y) 
\label{naturally reductive}
\end{equation}
Equivalently, the torsion $\widehat T$ of the normal connection is totally skewsymmetric.
\end{defi}

For a 3-symmetric space, the intrinsic connection coincides with the normal connection. As a consequence,
\begin{prop}
A Riemannian 3-symmetric space is nearly Kähler if and only if it is naturally reductive.
\end{prop}

Condition (\ref{naturally reductive}) is satisfied when $g$ is the restriction of an ${\rm Ad}(G)$-invariant scalar product $q$ on $\g$, representing a biinvariant metric on $G$. For example, on a compact, semi-simple Lie group $G$, the Killing form $B$ is negative definite, so we can take $q=-B$.

As we shall see in more detail, 3-symmetric spaces are related to twistor geometry :
\begin{prop}
The twistor space $M$ of a quaternion-Kähler (or Einstein, self-dual, in dimension 4) symmetric space is a 3-symmetric space with canonical almost complex structure $J_2$ (see theorem \ref{QK twistor is NK}).
\label{sym twistor is 3sym}
\end{prop}

This is compatible with the fact \cite{al,ee} that $(M,g,J_2)$ is quasi-Kähler, for a general twistor metric $g$.

\section{Twistor spaces and holonomy}

Francis Burstall \cite{bur2} practically resolved the problem of finding twistor spaces (complex manifolds fibring over a real manifold, such that the fibres are complex submanifolds) over Riemannian manifolds using the holonomy (for the moment, we are talking about the {\it Riemannian} holonomy). Here is a summary of his results that will motivate our own hypothesis. 

For an oriented Riemannian manifold $(M,g)$ of dimension $2n$, let $\pi : \ZZ \to M$ be the bundle whose sections are the almost complex structures on $M$, compatible with the metric $g$ and the orientation. It can be seen as the associated bundle of the metric structure with fibre 
\[ \mathcal Z(n) = \{ {\sf j} \in SO(2n) \ | \ {\sf j}^2 = -Id \} \simeq SO(2n)/U(n) \]
The manifold $\ZZ$ is itself equipped with an almost complex structure, defined in the following way. First, the fibres have a natural complex structure $\mathcal J^v$, because $SO(2n)/U(n)$ is an Hermitian symmetric space. Then, the Levi-Civita connection of $M$ induces a horizontal distribution $\mathcal H$ on $\ZZ$ that allows us to complete $\mathcal J^v$ into an almost complex structure $\mathcal J$ on the whole $T \ZZ$, setting
\[ \forall j \in \ZZ, \quad \J_j|_{\V} = \J^v_j, \quad \J_j|_{\H} = \pi^*j \]
where the point $j \in \ZZ$ is viewed as a complex structure on $T_xM$, $x=\pi(j)$. It is a celebrated result \cite{at} that $\J$ is integrable for a large class of 4-dimensional manifolds, called the "self-dual" manifolds. Unfortunately, in higher dimensions $2n \geq 6$, the Riemannian manifolds $M$ for which $\ZZ$ is a complex manifold are much rarer. In fact they must be locally conformally flat. We then look for a twistor space on $M$ among the {\it submanifolds} of $\ZZ$ i.e. we look for an almost complex submanifold $\Z$ of $\ZZ$ where $\J$ is integrable. To be an almost complex submanifold, $\Z$ must verify two conditions : \\
(\romannumeral 1) the fibre $\Z_x$ must be a complex submanifold of $\ZZ_x$. \\
(\romannumeral 2) the restriction to $\Z$ of the horizontal distribution $\H$ associated to the Levi-Civita connection $\nabla$ must be tangent to $\Z$. \\
For condition (\romannumeral 2) to be satisfied, it suffices to take a subbundle associated to the holonomy reduction. Consequently, if $H \subset SO(2n)$ is the holonomy group of $\nabla$, a natural idea is to consider the bundle with typical fibre the orbit of ${\sf j} \in \mathcal Z(n)$ under $H$. We then get a well behaved submanifold and condition (\romannumeral 1) is equivalent to an algebraic condition involving $H$. Secondly, a point $j \in \ZZ$ is in the zero set of the Nijenhuis tensor, $\N$, of $\J$ if and only if
\begin{equation} 
R_x(T^+_j,T^+_j)T^+_j \subset T^+_j
\label{J integrable}
\end{equation}
where $x=\pi(j)$ and $T^+_j \subset T_x^{\CM}M$ is the set of (1,0) vectors at $x$ with respect to $j$. Of course, according to the Newlander-Nirenberg theorem, $\J$ is integrable on $\Z$ if and only if (\ref{J integrable}) is satisfied for all $j \in \Z$.

Starting from Berger's list of holonomy representations of irreducible Riemannian manifolds \cite{berg}, F. Burstall investigates all cases. The case of symmetric spaces is very interesting and studied in details in the book \cite{bur} and the article \cite{bur3}. Indeed, symmetric spaces often admit several {\it flag manifolds} as twistor spaces. Concerning the non locally symmetric spaces, only the case where $H=Sp(q)Sp(1)$ ($n=2q$) gives fully satisfactory results : all quaternion-Kähler manifolds admit a complex twistor space with fibre $\CM P(1)$. In fact, they can be viewed as an analog, in dimension $4q$, $q \geq 2$, of Einstein, self-dual manifolds in dimension 4. For $U(n)$, i.e. for Kählerian manifolds, there is still a large class of manifolds having a twistor space with fibre $\CM P(n-1)$, the Bochner-flat Kähler manifolds (see \cite{ob}). Finally, the method gives no new results for $H$ being $Spin(7)$, $SU(n)$ or $Sp(q)$ (the manifold has to be conformally flat in order to admit an integrable twistor space). Of course the representation of $G_2$ is odd (seven) dimensional so it is not to be considered.

Consequently, we shall now interest ourselves to $H$-structures {\it with torsion}. But in this case $\nabla$ is not the appropriate connection anymore since it is not a $H$-connection. Let $\nt$ be a $H$-connection. We construct a new almost complex structure $\Jt$ on $\ZZ$ using the horizontal distribution of $\nt$ and a submanifold $\Z$ of $\ZZ$, as before :
\begin{defi}
Let $M$ be a $H$-manifold of dimension $2n$, ${\sf j} \in \mathcal Z(n)$. The $H$-twistor space associated to ${\sf j}$ is the associated bundle $\Z$ of the $H$-structure with fibre the $H$-orbit $\mathcal O$ of ${\sf j}$ in $\mathcal Z(n)$.
\end{defi}

\begin{NB} 
We shall talk also of the $H$-twistor space associated to $\mathcal O$, or to a point $j \in \ZZ$, identifying $T_x M$ to $\RM^{2n}$ -- and thus $\ZZ_x$ to $\mathcal Z(n)$ -- by means of a frame $p \in H_x(M)$. This is independent on the choice of $p$, because we consider $H$-orbits.
\end{NB}

Since $\Z$ is an associated bundle of the $H$-structure, the horizontal distribution of $\nt$ on $\ZZ$ is tangent to it. Then, $\Z$ is an almost complex submanifold of $(\ZZ, \Jt)$ if and only if the vertical distribution is also stable by $\Jt$, i.e. $\mathcal O$ is a complex submanifold of $\mathcal Z(n)$.

\begin{prop}
Let ${\sf j} \in \mathcal Z(n)$. Let $H$ be a subgroup of $SO(2n)$ acting on $\mathcal Z(n)$ with stabilizer $K$ at ${\sf j}$. Then the orbit of ${\sf j}$, isomorphic to $H/K$, is a complex submanifold of $\mathcal Z(n)$ if and only if ${\sf j}$ belongs to $N(H)$, the normalizer of $H$ in $SO(2n)$.
\label{stabilizer of H}
\end{prop}

\begin{proof}
We can always assume $U(n)$ is the stabilizer of our fixed ${\sf j}$. Let $\mathfrak{so}(2n)$, $\mathfrak{u}(n)$ be the Lie algebra of $SO(2n)$, $U(n)$, respectively.
\[ \mathfrak{so}(2n) = \mathfrak{u}(n) \oplus \mathfrak z \]
where $\mathfrak z$ represents the tangent space of $\mathcal Z(n) \simeq SO(2n)/U(n)$. Moreover, $\mathfrak{so}(2n)$ may be identified with the space of 2-forms on $\RM^{2n}$ and if we complexify :
\[ \mathfrak{so}(2n)^{\CM} = \mathfrak{u}(n)^{\CM} \oplus \mathfrak z^{2,0} \oplus \mathfrak z^{0,2} \]
where $\mathfrak{u}(n)^{\CM}$ is identified with the space of {\it complex} (1,1)-forms on $\RM^{2n}$ and $\mathfrak z^{2,0}$ (resp. $\mathfrak z^{0,2}$) is identified with the space of 2-forms of type (2,0) (resp. (0,2)). This is also the eigenspace decomposition, for the eigenvalues $0$, $2i$, $-2i$, of ${\sf j}$, acting as an element of the Lie algebra. The natural complex structure on $\mathcal Z(n)$ is defined to be the invariant almost complex structure on the homogeneous space represented by $\mathfrak z^{2,0} \subset \mathfrak z^{\CM}$. It is integrable because ${\sf j}$ acts on $[\mathfrak z^{2,0},\mathfrak z^{2,0}]$ with eigenvalue $2i+2i=4i$ so this space can only be $\{0\}$. Now, if ${\sf j} \in N(H)$, it preserves $H$ or $\mathfrak h$, so there is a decomposition into eigenspaces :
\begin{equation} 
\h^{\CM} = \k^{\CM} \oplus \v^{2,0} \oplus \v^{0,2} 
\label{H/K complex submanifold} 
\end{equation}
where $\k^{\CM} = \mathfrak{u}(n)^{\CM} \cap \h^{\CM}$, $\v^{2,0} = \mathfrak z^{2,0} \cap \h^{\CM}$, $\v^{2,0} = \mathfrak z^{2,0} \cap \h^{\CM}$ ($K \subset U(n)$ because it is the stabilizer of ${\sf j}$ in $H$). Consequently $H/K$ is a complex submanifold of $SO(2n)/U(n)$.

Conversely, if $H/K$ is a complex submanifold of $SO(2n)/U(n)$, (\ref{H/K complex submanifold}) holds so ${\sf j}$ preserves $\h$, i.e. belongs to $N(H)$.
\end{proof}

The following corollary will be of considerable importance in the sequel.
\begin{coro}
Assume $H$ is a semi-simple subgroup of $SO(2n)$. If $\mathcal O \simeq H/K$ is a complex submanifold of $\mathcal Z(n)$, then $K$ is the centralizer of a torus in $H$.
\label{K centralizer in H}
\end{coro}

Equivalently, $K$ preserves an element ${\sf j}_0 \in \h$ or $\mathcal O$ is isomorphic to an adjoint orbit of $H$.

\begin{proof}
Let $C(H)$ be the centralizer of $H$ in $SO(2n)$ :
\[ C(H) = \{ g \in SO(2n) \ | \ \forall h \in H, ghg^{-1} = h \} \]
Denote by $\mathfrak{c}(\h)$ its Lie algebra and by $\n(\h)$ the Lie algebra of $N(H)$, we have
\begin{equation}
\n(\h) = [\h,\h] \oplus \mathfrak{c}(\h) 
\label{n(h) = h + c(h)}
\end{equation}
The proof of this fact is given in \cite{re}, p15 and uses an ${\rm Ad}(H)$-invariant scalar product on $\n(\h)$.

Now, if $H/K$ is a complex submanifold of $SO(2n)/U(n)$, by proposition \ref{stabilizer of H}, ${\sf j}$ belongs to $N(H)$, or $\n(\h)$. Thus, it can be decomposed according to (\ref{n(h) = h + c(h)}) :
\[ {\sf j} = {\sf j}_0 + {\sf j}_1 \]
where ${\sf j}_0 \in \h$ and ${\sf j}_1 \in \mathfrak{c}(\h)$. Moreover an element $h \in H$ stabilizes ${\sf j}$ if and only if it stablizes ${\sf j}_0$ because it always stabilizes ${\sf j}_1$ by the definition of the centralizer. Consequently $K$ is the centralizer of ${\sf j}_0 \in H \cap \h$, or of the torus generated by ${\sf j}_0$.
\end{proof}

In this article we are mainly concerned with the case where $H$ is a subgroup of $U(n)$. Then $M$ is an almost Hermitian manifold $(M,g,J)$ and the necessary condition \ref{stabilizer of H} becomes :
\begin{coro}
Let $H$ be a subgroup of $U(n)$ such that $N(H) \subset U(n)$. Let $\Z$ be the $H$-twistor space over a $H$-manifold $M$ associated to ${\sf j} \in \mathcal Z(n)$. Suppose that $\mathcal O = H.{\sf j}$ is a complex submanifold of $\mathcal Z(n)$ and so $\Z$ is an almost complex submanifold of $(\ZZ,\Jt)$, where $\Jt$ is the almost complex structure associated to a $H$-connection $\nt$. Then, $\forall j \in \Z$, $j$ commutes with $J_x$ on $T_x M$, $x=\pi(j)$.
\label{j commute avec J0}
\end{coro}

\begin{proof}
It is important to understand that in the hypothesis $H \subset U(n)$ of corollary \ref{j commute avec J0}, $U(n)$ is not the stabilizer of ${\sf j}$ anymore (otherwise $\Z$ would be trivial) but of another ${\sf j}_0 \in \mathcal Z(n)$. Now, if $H/K$ is a complex submanifold, ${\sf j} \in N(H) \subset U({\sf j}_0)$ so it commutes with every element of the center and in particular with ${\sf j}_0$ itself. In the geometrical background, this translates exactly into corollary \ref{j commute avec J0}.
\end{proof}

\begin{NB}[Structure of the reduced twistor space]
An important example is $H=U(n)$ itself. Then $N(H)=U(n)$ and $H/K$ being a complex submanifold is {\sl equivalent} to ${\sf j}$ commuting with ${\sf j}_0$. This is why we considered, in \cite{bur2}, the "reduced twistor space" ${\boldsymbol{{\boldsymbol{\mathcal Y}}}} \subset \ZZ$ of an almost complex manifold $(M,g,J)$ whose fibre at $x \in M$ is exactly the set of complex structures $j : T_xM \to T_xM$ commuting with $J_x$. More precisely, let $j \in {\boldsymbol{\mathcal Y}}$, $x=\pi(j)$. We diagonalize $j$ and $J_x$ simultaneously i.e. we find 4 subspaces $R_j^+$, $R_j^-$, $S_j^+$ and $S_j^-$ of $T^{\CM}_xM$ such that the eigenspaces of $j$ for the respective eigenvalues $i$, $-i$ are
\[ T^+_j = R^+_j \oplus S^+_j, \quad T^-_j = R^+_j \oplus S^+_j \]
and the eigenspaces of $J_x$ :
\[ T^+_xM = R^+_j \oplus S^-_j, \quad T^-_xM = R^-_j \oplus S^+_j \]
Thus $R^-_j = \overline{R^+_j}$ and $S^-_j = \overline{S^+_j}$. Equivalently, if we denote
\[ R_j=T_xM \cap (R^+_j \oplus R^-_j) \quad \text{and} \quad  S_j=T_xM \cap (S^+_j \oplus S^-_j), \]
the endomorphisms $j$ and $J_x$ coincide on $R_j$ but have opposite signs on $S_j$. Then $j$ is determined by $R^+_j$ or $R_j$ and the connected components of ${\boldsymbol{{\boldsymbol{\mathcal Y}}}}$ are characterized by the (complex) dimension of $R^+$. Denote
\begin{equation}
{\boldsymbol{\mathcal Y}}_r = \{j \in \ZZ \ | \ j \circ J_x = J_x \circ j, \quad \mathrm{dim}_{\CM} R^+_j = r \}, \quad r = 1, \ldots n-1
\label{Grassmannian bundles} 
\end{equation}
The group $U(n)$ acts transitively on each fibre of ${\boldsymbol{\mathcal Y}}_r$ which is thus an $U(n)$-almost complex twistor space over $M$, with fibre the complex Grassmannian $Gr_r(\CM^n) \simeq U(n)/U(r) \times U(n-r)$. 
\end{NB}

Now, let $\Nt$ be the Nijenhuis tensor of $\Jt$, the conditions for the vanishing of $\Nt_j$ are more complicated in this case. In fact (\ref{J integrable}) is only a particular case of a more general system established by O'Brian, Rawnsley \cite{ob} :
\begin{prop}[O'Brian \& Rawnsley]
Let $M$ be a Riemannian manifold and $\Jt$ the almost complex structure on $\ZZ$ associated to a metric connection $\nt$. Let $\widetilde T$, $\Rt$ be the torsion and curvature of $\nt$. A point $j \in \ZZ_x$ lies in the zero set of the Nijenhuis tensor $\Nt$ of $\Jt$ if and only if
\begin{equation} 
\widetilde T_x(T^+_j,T^+_j) \subset T^+_j
\label{Jt integrable torsion}
\end{equation}
\begin{equation} 
\Rt_x(T^+_j,T^+_j)T^+_j \subset T^+_j
\label{Jt integrable curvature}
\end{equation}
where $T^+_j \subset T_x^{\CM} M$ is the set of (1,0) vectors with respect to $j : T_x M \to T_x M$.
\end{prop}

\section{Fibration of degree 3}

Let $M$, in this section, be a homogeneous space $G/H$, where $G$ is a connected Lie group and $H$ is a closed subgroup of $G$. Assume that $G/H$ is reductive and let $g$ be an ${\rm Ad}(H)$-invariant scalar product on $\m$ (so that $H \subset SO(n)$), inducing an invariant metric on $M$ denoted the same way. Just like the theory of twistor spaces of {\it inner} symmetric spaces is simpler than in the outer case, we shall make here the simplification that the groups $G$ and $H$ have same rank. We also assume that $G$ is compact, semi-simple. 

The theory of twistor spaces of "inner" homogeneous spaces relates to the study of flag manifolds :

\begin{defi}
A flag manifold is a homogeneous space $G/K$ where $K$ is the centralizer of a torus in $G$.
\end{defi}

\begin{prop}
Let $M=G/H$ be a homogeneous space such that $\mathrm{rank} \ G = \mathrm{rank} \ H$. Let $\Z$ be the $H$-twistor space associated to ${\sf j} \in \mathcal Z(n)$. Suppose that $\Z$ is an almost complex manifold of $(\ZZ,\Jt)$. Then $\Z$ is a flag manifold.
\label{Z is flag}
\end{prop}

\begin{proof}
For a homogeneous space, the $H$-structure is identified to $G$ itself, and the orbit of ${\sf j}$ is isomorphic to $H/K$, where $K$ is the stabilizer of ${\sf j}$ in $H$, so
\[ \Z = G \times_H H/K \simeq G/K \]
Now, if $H/K$ is a complex submanifold of $\mathcal Z(n)$, $K$ is the centralizer of a torus in $H$ by corollary \ref{K centralizer in H}. In particular, it is a subgroup of maximal rank of $H$. But $G$ and $H$ have same rank so $K$ is a subgroup of maximal rank of $G$. Then, by the work of Wolf \cite{wo3}, 8.10, it is also the centralizer of a torus {\it in $G$}.
\end{proof}

\begin{NB}
In this homogeneous context we will talk of $H$-twistor spaces associated to $\n^+$, instead of ${\sf j} \in \mathcal Z(n)$, where $\n^+$ is a maximal isotropic subspace of $\m^{\CM}$, representing a complex structure on $\m$ compatible with $g$. Indeed, we see an element $\varphi$ of $G$ as a frame $\varphi : \m \simeq T_o M \to T_{\pi(g)} M$, identifying $\varphi$ with the differential of the left action $[\varphi'] \mapsto [\varphi \varphi']$ on $M$. Thus, to work directly with $G$ as a $H$-structure on $M$, it is more appropriate to replace $\RM^{2n}$ by $\m$ in the definition of $\mathcal Z(n)$.

In other words, let $j$ be a complex structure on $T_x M$, compatible with $g_x$. We transport $j$ or $T^+_j$ by $\varphi \in G$, where $x=[\varphi]$, to obtain a complex structure $\n^+$ on $T_o M \simeq \m$. Another choice of representant of $x$, $\varphi'=\varphi h$, $h \in H$ would give a different complex structure $\n'^+ \subset \m^{\CM}$ but related to $\n^+$ via
\[ \n'^+ = {\rm Ad}_h \n^+. \]
So to be precise we should refer to the $H$-twistor space associated to an {\it orbit} ${\rm Ad}(H)\n^+$ of complex structures on $\m$.
\end{NB}

Moreover we can choose for the $H$-connection $\widehat \nabla$, the normal connection of the reductive homogeneous space. We denote by $\widehat J$ the almost complex structure on $\ZZ$ associated to it. Then, since the torsion and curvature of $\widehat \nabla$ are $G$-invariant, integrability conditions for $\widehat J$ reduce to conditions on $\n^+$ :

\begin{prop}
Let $M \simeq G/H$ be a reductive homogeneous space. Let $\Z$ be the $H$-twistor space on $M$ associated to $\n^+$. The following conditions are equivalent : \\
(\romannumeral 1) $\widehat J$ is integrable on $\Z$. \\
(\romannumeral 2) $\exists j \in \Z$ such that $\widehat \N_j = 0$, where $\widehat \N_j$ is the Nijenhuis tensor of $\widehat J$. \\
(\romannumeral 3) The subset $\n^+$ satisfies
\begin{equation}
[\n^+,\n^+]^{\m} \subset \n^+, \quad [\,[\n^+,\n^+]^{\h},\n^+] \subset \n^+
\label{Jt integrable homogene}
\end{equation}
where superscript $\m$, $\h$ stand for the projections on the appropriate subspace.
\label{twistor homogeneous}
\end{prop}

\begin{proof}
As we said, equivalence between (\romannumeral 1) and (\romannumeral 2) comes from the invariance of $\widehat T$, $\widehat R$. As for (\ref{Jt integrable homogene}), this is nothing else than (\ref{Jt integrable torsion}), (\ref{Jt integrable curvature}), using the expressions of $\widehat T$, $\widehat R$ given in lemma \ref{normal torsion curvature}.
\end{proof}

\vs

Thus, we will now interest ourselves to flag manifolds. Let $T$ be a maximal torus of $G$ and $S$ a subtorus of $T$. Denote by $\t$, $\s$ their Lie algebras and by $\mathcal B$ a choice of simple roots for the root system $\mathcal R$ associated to $T$. 

Let $\mathcal A$ be the subset of $\mathcal B$ whose complement is
\[ \overline{\mathcal A} = \{ \alpha \in \mathcal B \ | \ \forall X \in \s, \ \alpha(X) = 0 \} \]

\begin{lemm}
The centralizer $K$ of $S$ is the subgroup of $G$ with Lie algebra $\k = \k_{\A}$.
\label{k=kA}
\end{lemm}
\begin{proof}
Recall the definition (\ref{kA,lA+}). Equivalently, $\k_{\A}^{\CM}$ is the sum of all root spaces $\g_{\alpha}$ where $\alpha$ is a sum of simple roots $\beta \in \overline{\mathcal A}$. Let $X \in T$, $Y \in \g_{\alpha}$. The commutator $[X,Y]$ is given by
\[ {\rm ad}(X)Y = 2\pi i \alpha(X) = 2\pi i \sum_{\beta \in \mathcal B} n_{\beta}(\alpha) \beta (X), \]
where all $n_{\beta}$ have same sign. It is zero for all $X \in S$, and thus $Y \in \k$ if and only if $n_{\beta} = 0$ whenever $\beta \in \mathcal A$.
\end{proof}

Conversely, for each subset $\A$ of $\mathcal B$, $\k_{\A}$ is the centralizer of $\s=\{X \in \t \ | \ \forall \alpha \in \overline{\A}, \ \alpha(X)=0 \}$.

Now, there is also a parabolic subalgebra $\p_{\A}$ (and a parabolic subgroup $P_{\A}$) associated to the subset $\A$ as in (\ref{pA}). Denote, as in the preliminaries, $\l_{\A}^+$ its nilradical, $\l_{\A}^- = \overline{\l_{\A}^+}$ and
\[ \l_{\A} = \g \cap (\l_{\A}^+ \oplus \l_{\A}^-) \]
so that $\l_{\A}^+$ is a complex structure on $\l_{\A}$,
\[ \g = \k_{\A} \oplus \l_{\A}, \quad [\k_{\A},\l_{\A}] \subset \l_{\A} \]
and
\[ [\k_{\A},\l_{\A}^+] \subset \l_{\A}^+.\]
Then $\l_{\A}$ is an ${\rm Ad}(K)$ complement of $\k$ in $\g$, $\l_{\A}^+$ represents an invariant almost complex structure on $G/H$. Moreover, by lemma \ref{lA+ nilradical}, $\l^+_{\A}$ is a subalgebra so this almost complex structure is in fact integrable, corresponding to the isomorphism
\begin{equation} 
G/K_{\A} \simeq G^{\CM}/P_{\A} 
\label{G/K=G^C/P}
\end{equation}

Conversely, let $\p$ be parabolic subalgebra corresponding to a parabolic subgroup $P$. Then, $G^{\CM}/P$ is a flag manifold since $\p$ might be written $\p_{\A}$, for a base $\mathcal B$ and a subset $\mathcal A \subset \mathcal B$ according to theorem \ref{parabolic subalgebras} and the decomposition (\ref{p=k^c + l-}) is still valid. 

Note that the isomorphism (\ref{G/K=G^C/P}) depends on the choice $\mathcal B$ of simple roots. It is only in this second form, $G^{\CM}/P$, that the flag manifold reveals as the twistor space of a symmetric space, as shown by Burstall and Rawnsley, but also of a 3-symmetric space, as we will now see.

We fix now such a complex structure -- or such a $\p$ -- on our flag manifold $G/K$ and denote $\p_{\A}$, $\l_{\A}$, $\l_{\A}^+$, $\l_{\A}^-$ simply by $\p$, $\l$, $\l^+$ and $\l^-$.

From the nilradical of $\p$, we define a finite series :
\begin{defi}
The first canonical series associated to $\p$ is the central descending series $(\l_i)_{i \geq 1}$ of $\l^+$ given by :
\[ 
\l_1 = \l^+, \quad \l_2 = [\l_1,\l_1], \quad \l_3 = [\l_1,\l_2], \ldots 
\]
Then, the second canonical series $(\g_i)_{i \geq 1}$ is defined by $\g_i = \l_i \cap \l_i^{\perp}$, where the orthogonal is taken with respect to the Killing form of $G$. Finally, we put $\g_0 = \k^{\CM}$ and for $i \leq 0$, $\g_i = \overline{\g_{-i}}$.
\end{defi}

First, by definition
\[ \g^{\CM} = \sum_{i \in \ZM} \g_i, \quad \p = \sum_{i \geq 0} \g_i, \quad \l^+ = \sum_{i > 0} \g_i \]
Then, the $\g_i$ have the following fundamental property,
\begin{theo}[Burstall, Rawnsley]
Let $G$ be a compact semi-simple Lie group, $\p$ a parabolic subalgebra of $\g$. The second canonical series associated to $\p$ satisfies
\[ \forall i,j \in \ZM, \quad [\g_i,\g_j] \subset \g_{i+j} \]
Moreover $\g_1$ generates the series in the sense that
\begin{equation}
\g_r = [\g_1,[\g_1,[\ldots,\g_1]\ldots]], 
\label{gr=[g1,[..]]}
\end{equation}
$\g_1$ appearing $r$ times in the last formula.
\label{canonical series}
\end{theo}

The proof is by exhibiting an $X \in \k$, that they called {\it the canonical element of $\p$}, such that $\g_k$ is the eigenspace for the eigenvalue $k.i$ of ${\rm ad}(X)$.

We will now construct another homogeneous space using the second canonical series. Let
\begin{equation}
\h = \sum_{i \in \NM} \leftr \g_{3i} \rightr, \quad \m^+ = \sum_{i \in \ZM} \g_{3i+1}, \quad \m^- = \sum_{i \in \ZM} \g_{3i+2}
\label{3-fibration}
\end{equation}
so that $\m^+ = \overline{\m^-}$, because $-(3i+1)=-3(i+1)+2$ and if we put $\m = \g \cap (\m^+ \oplus \m^-)$, $\m^+$ represents a complex structure on $\m$. 

First, $\h$ is a subalgebra because $3i + 3j = 3(i+j)$ for all $i,j \in \ZM$ so $[\g_{3i},\g_{3j}] \subset \g_{3(i+j)} \subset \h$ by theorem \ref{canonical series}. Let $H$ be the subgroup of $G$ corresponding to $\h$. In particular, $\k = \leftr g_0 \rightr$ is contained in $\h$ so $K$ is a subgroup of $H$.

\begin{defi}
The fibration $G^{\CM}/P \simeq G/K \to G/H$ is called fibration of degree 3 of the (complex) flag manifold.
\end{defi}

Furthermore,
\[ [ \h,\m^+ ] \subset \m^+ \quad \text{and} \quad [ \h,\m^- ] \subset \m^- 
\]
because $3i + (3j+1) = 3(i+j) + 1$ and $3i + (3j+2) = 3(i+j)+2$, for all $i,j \in \ZM$, and as a consequence,
\[ [\h,\m] \subset \m. \]
We also have
\[ \h^{\CM} = \sum_{i \in \ZM} \g_{2i}, \quad \g^{\CM} = \h^{\CM} \oplus \m^+ \oplus \m^-, \quad \g = \h \oplus \m \]
so $\m$ is an ${\rm Ad}(H)$-complement of $\h$ in $\g$ and $\m^+$ defines an invariant almost complex structure on the reductive homogenous space $G/H$.

\begin{prop}
The homogeneous space $G/H$ is a 3-symmetric space with canonical almost complex structure represented by $\m^+ \subset \m^{\CM}$.
\end{prop}
\begin{proof}
We easily verify the equations (\ref{3symmetric}), using theorem \ref{canonical series}.
\end{proof}

\begin{NB} 
This is very similar to the construction of the canonical fibration by Burstall, Rawnsley, going from a flag manifold to a {\it symmetric} space. We should think of the canonical fibration as a "fibration of degree 2" : let $\mathfrak f = \sum_{i \in \NM} \leftr \g_{2i} \rightr$ and $\mathfrak q = \sum_{i \in \NM} \leftr \g_{2i+1} \rightr$. By theorem \ref{canonical series}, $\mathfrak f$ is a subalgebra, corresponding to a subgroup $F$ of $G$ ; $\mathfrak q$ is an ${\rm Ad}(H)$-invariant complement of $\mathfrak f$ in $\g$ and $[\mathfrak q,\mathfrak q] \subset \mathfrak f$ ($\forall i,j \in \ZM$, $(2i+1) + (2j+1)=2(i+j+1)$) so $G/F$ is a symmetric space.
\end{NB}

Since $G$ is compact, semi-simple, the Killing form $B$ is negative definite so the restriction of $-B$ to $\m$ is an ${\rm Ad}(H)$-invariant scalar product that makes $M$ a naturally reductive, strictly nearly Kähler, homogeneous space.

Moreover, $\p$ gives a complex structure on $\m$ :
\begin{equation}
\n^+ = \p \cap \m^{\CM} = \bo_{i \geq 0} \g_{3i+1} \oplus \bo_{i \geq 0} \g_{3i+2}
\label{n+}
\end{equation}

\begin{prop}
The space $G^{\CM}/P$ is the $H$-twistor space on $G/H$ associated to $\n^+$, with integrable complex structure associated to the intrinsic connection.
\label{the 3-fibration is twistor}
\end{prop}

\begin{proof}
The fibre $H/K$ of the fibration of degree 3 is a complex submanifold of $G^{\CM}/P$. Indeed, an ${\rm Ad}(K)$-complement of $\k$ in $\h$ is $\mathfrak a = \sum_{i > 0} \leftr \g_{2i} \rightr$ and the subspace $\sum_{i > 0} \g_{2i}$ of $\mathfrak a^{\CM}$ is a subalgebra of $\l^+$, so it defines a complex structure on $H/K$ which coïncides with the restriction of $\l^+$. Consequently $G/K \to G/H$ is a twistor fibration.
\end{proof}

We will soon prove a partial converse of proposition \ref{the 3-fibration is twistor}. Even if our main goal are the nearly Kähler 3-symmetric spaces, the first results are more general. We shall assume that $M \simeq G/H$ is a Riemannian, reductive homogeneous space such that $\m$ is the orthogonal of $\h$ for $q=-B$, where $B$ is the Killing form of $G$, and $g$ equals the restriction to $\m$ of $q$. As a consequence, $M$ is naturally reductive. We also assume as before that $G$ and $H$ have same rank. Then, if $\Z \simeq G/K$ is an almost complex manifold, it is a flag manifold by proposition (\ref{Z is flag}). But if $\Z$ is a {\it complex} manifold we have more : a complex structure on $G/K$, i.e. an identification of $G/K$ with $G^{\CM}/P$, for some parabolic subgroup $P$ of $G^{\CM}$.

\begin{theo}
Let $\n^+$ be a maximal isotropic subspace of $\m^{\CM}$. Then $\n^+$ satisfies (\ref{Jt integrable homogene}) if and only if $[\n^+,\n^+] + \n^+$ is the nilradical of a parabolic subalgebra $\p \subset \g^{\CM}$. As a consequence, the twistor space associated to $\n^+$ is integrable, isomorphic to $G^{\CM}/P$ where $P$ denotes the parabolic subgroup of $G^{\CM}$ corresponding to $\p$. 
\label{G^c/P -> G/H}
\end{theo}

\begin{proof}
We were very much inspired by the proofs of Theorem 4.8, p44 and Lemma 5.1, p64 of \cite{bur}. 

Let 
\[ \l^+ = [\n^+,\n^+] + \n^+ = [\n^+,\n^+]^{\h} \oplus \n^+ \]
where the second decomposition respects the sum $\g=\h \oplus \m$. Firstly, $[\n^+,\n^+] \subset [\l^+,\l^+]$ is contained in $\l^+$ (which has thus a chance to be a subalgebra) if and only $[\n^+,\n^+]^{\m} \subset \n^+$. This corresponds to the first part of (\ref{Jt integrable homogene}). Secondly, $[\h,\m] \subset \m$ so $[[\n^+,\n^+]^{\h},\n^+] \subset \l^+$ if and only if it is contained in $\n^+$ and this corresponds to the second part of (\ref{Jt integrable homogene}). Moreover, this is equivalent to $[[\n^+,\n^+],\n^+] \subset \l^+$ provided that $[\n^+,\n^+]^{\m} \subset \n^+$. It remains to show that (\ref{Jt integrable homogene}) implies
\[ [[\n^+,\n^+],[\n^+,\n^+]] \subset [\n^+,\n^+] + \n^+ \]
By the Jacobi identity,
\begin{eqnarray*}
[[\n^+,\n^+],[\n^+,\n^+]] & \subset & [[[\n^+,\n^+],\n^+],\n^+] \\
& \subset & [\n^+,\n^+] + [[\n^+,\n^+],\n^+] \\
& \subset & [\n^+,\n^+] + \n^+
\end{eqnarray*}
Finally, $\l^+$ is a subalgebra if and only if $\n^+$ satisfies (\ref{Jt integrable homogene}).

Let $K$ be the stabilizer of $\n^+$ and let $\k$ be its Lie algebra. An element $X \in \h$ belongs to $\k$ if and only if
\[ [X,\n^+] \subset \n^+ \]
Since $[\h,\m] \subset \m$ and $\n^+$ is a maximal isotropic subset of $\m^{\CM}$ with respect to $g$, this is equivalent to
\[ g([X,\n^+],\n^+) = 0 \]
But $g$ is the restriction of $q$, which is ${\rm Ad}(G)$-invariant, so $X \in \k$ if and only if
\[ q(X,[\n^+,\n^+]^{\h}) = 0 \]
Consequently, denoting $\l^- = \overline{\l^+}$, we have
\[ \g^{\CM}= \k^{\CM} \oplus \l^+ \oplus \l^- \]

We want to prove that
\[ \p = \k^{\CM} \oplus \l^+ \]
is a parabolic subalgebra.  To do this, we need to find a base $\mathcal B$ of a root system, and a subset $\A$ of $\mathcal B$ such that $\p = \p_{\A}$. 

Before that, we prove that $\p$ is a subalgebra. By the Jacobi identity : 
\[ [\k^{\CM},[\n^+,\n^+]] \subset [[\k^{\CM},\n^+],\n^+] \subset [\n^+,\n^+] \subset \l^+ \]
We also prove that $\l^+$ is isotropic. Let $\v^+ = [\n^+,\n^+]^{\h}$. Since $\n^+$ is maximal isotropic, the second part of (\ref{Jt integrable homogene}) is equivalent to
\[ g([\v^+,\n^+],\n^+) = 0 \]
But $g$ is the restriction of an ${\rm Ad}(G)$-invariant scalar product $q$, and $\h^{\CM}$, $\m^{\CM}$ are $q$-orthogonal so this is equivalent to
\[ q(\v^+,\v^+)=0. \]
The last equation tells exactly that $\v^+$ is isotropic and so, so is $\l^+ = \v^+ \oplus \n^+$. 

By proposition \ref{Z is flag}, if $\Z$ is an almost complex manifold, $K$ is the centralizer of a torus. Thus, it contains a maximal torus $T \subset G$. We consider the root system $\mathcal R$ associated to $T$. We define a subset $\mathcal S$ of $\mathcal R$ by
\[ \l^+ = \sum_{\alpha \in \mathcal S} \g_{\alpha} \]
This subset has the following two properties : \\
(\romannumeral 1) $\mathcal S \cap - \mathcal S = \{0 \}$. \\ Indeed $\l^+$ is isotropic so we can't find $\alpha \in \mathcal R$ such that $\g_{\alpha} \subset \l^+$ {\it and} $\overline{\g_{\alpha}} \subset \l^+$. And : \\
(\romannumeral 2) If $\alpha$, $\beta$ belong to $\mathcal S$, then $\alpha + \beta \notin \mathcal R$ or $\alpha + \beta \in \mathcal S$ ($\mathcal S$ is {\it closed}). \\ This comes from $\l^+$ being a subalgebra.

By a result of \cite{bo} cited in \cite{bur}, p27, if a roots subset satisfies (\romannumeral 1) and (\romannumeral 2) then we can find a base $\mathcal B$ such that $\mathcal S \subset \mathcal R^+$, the set of positive roots. With this choice of simple roots, there exists $\A \subset \mathcal B$ such that $\k = \k_{\A}$ (cf lemma \ref{k=kA}). Of course, if $\alpha \in \mathcal S$, $n_{\A}(\alpha) \geq 0$ because it is a positive root so
\[ \l^+ \subset \sum_{n_{\A}(\alpha) \geq 0} \g_{\alpha} \]
and finally we must have $\l^+=\l_{\A}^+$.
\end{proof}

\begin{NB} It would be enough to assume that $g$ is the restriction of {\it any} ${\rm Ad}(G)$-invariant scalar product on $\g$, i.e. $G/H$ is {\it normal}. However, this is equivalent to our assumption as soon as $G$ is simple.
\end{NB}

\begin{NB}
The fibration we obtain is trivial if and only if $[\n^+,\n^+]^{\h} = \{0\}$ and $[\n^+,\n^+] \subset \n^+$ i.e. $\n^+$ is already the nilradical of a subalgebra and $G/H$ is a flag manifold.
\end{NB}

The demonstration of theorem \ref{G^c/P -> G/H} gives us an explicit way to find $\n^+$ satisfying (\ref{Jt integrable homogene}) and so construct a complex twistor space over $G/H$. Indeed, we needed to show that $\l^+$, and so $\n^+$, are included in the sum of all positive root weight spaces, for a given base $\mathcal B$. Conversely, we have :

\begin{prop}
Let $G/H$ be a homogeneous space such that $\mathrm{rank} \ H = \mathrm{rank} \ G$. Let $T \subset H \subset G$ be a maximal torus. For a choice of simple roots $\mathcal B$, denote by $\n^+$ the complex structure on $\m$ given by
\begin{equation} 
\n^+ = \m^{\CM} \cap \bo_{\alpha \in \mathcal R^+} \g_{\alpha} 
\label{n+ associated to B}
\end{equation}
Then $\n^+$ satisfies (\ref{Jt integrable homogene}), i.e. the $H$-twistor space associated to $\n^+$ is complex.
\label{twistor space associated to B}
\end{prop}

\begin{proof}
The subset $\l^+_{\mathcal B} = \bo_{\alpha \in \mathcal R^+} \g_{\alpha}$ is a subalgebra of $\g^{\CM}$ so $[\n^+,\n^+] \subset \l^+_{\mathcal B}$, $[\n^+,\n^+]^{\m} \subset \n^+$ and $[[\n^+,\n^+]^{\h},\n^+] \subset \n^+$, by the definition of $\n^+$.
\end{proof}

\begin{NB}
In particular, suppose that $M=G/H$ is a 3-symmetric space with canonical complex structure given by $\m^+$. If there exists a fibration of degree 3 with base $M$, then this is one of the fibrations described in theorem \ref{G^c/P -> G/H} for the subset $\n^+$ of $\m^{\CM}$ given by (\ref{n+}). Indeed by (\ref{gr=[g1,[..]]}), $\forall i \geq 0$, $[\g^1,\g^{3i+2}] = \g^{3(i+1)}$ (with equality and not just inclusion) so 
\[ [\n^+,\n^+]^{\h} = \bo_{i \geq 1} \g_{3i}, \quad [\n^+,\n^+]^{\h} \oplus \n^+ = \p \]
However a fibration of degree 3 has the additional property, for the second canonical series associated to $\p$,
\begin{equation}
\g^1 \subset \m^+.
\label{fibration degree 3 g1,m+}
\end{equation}
Note that, if $G^{\CM}/P \to G/H$ is one of the fibration described in theorem \ref{G^c/P -> G/H}, we always have $\g^1 \subset \n^+$ because $\l^1 = \l^+ = \n^+ + [\n^+,\n^+]$ and so $\l^2$ contains $[\n^+,\n^+]$.
\label{g1,m+}
\end{NB}

\begin{prop}
Let $M \simeq G/H$ be a strictly nearly Kähler 3-symmetric space such that $G$ is semi-simple and the naturally reductive metric on $M$ is given by the restriction of $-B$. Let $G^{\CM}/P$ be a complex $H$-twistor space over $M$. Suppose moreover that the second canonical series associated to $\p$ and the canonical almost complex structure $\m^+$ satisfy (\ref{fibration degree 3 g1,m+}). Then, $G^{\CM}/P \to G/H$ is the fibration of degree 3 associated to $\p$.
\label{reciproque}
\end{prop}

\begin{proof}
Once we have $\g^1 \subset \m^+$, we easily get, by (\ref{3symmetric}),
\[ \g^2=[\g^1,\g^1] \subset [\m^+,\m^+] \subset \m^-, \ \g^3 = [\g^1,\g^2] \subset [\m^+,\m^-] \subset \h^{\CM}, \ \g^4 = [\g^1,\g^3] \subset \m^+, \ldots \]
So $\h,\m^+,\m^-$ are like in the definition (\ref{3-fibration}) of the fibration of degree 3.
\end{proof}

Finally we can replace condition (\ref{fibration degree 3 g1,m+}) by another, perhaps easier to apply : 
\begin{lemm}
Let $\mathcal B$ be a base of simple roots such that $\p = \p_{\A} \subset \bo_{\alpha \in \mathcal R^+} \g_{\alpha}$, $\A \subset \mathcal B$. Then $\g^1 \subset \m^+$ if and only if
\[ \forall \beta \in \A, \quad \g_{\beta} \subset \m^+ \]
\label{reciproque bis}
\end{lemm}
\begin{proof}
To show this, we first make the following simple remark that, since $\l^1 = \l^+_{\A}$ is given by (\ref{kA,lA+}),
the roots $\alpha$ whose weight spaces compose $\l^2=[\l^1,\l^1]$ all have $n_{\A}(\alpha) \geq 2$, etc., and so, for all $r \in \ZM$,
\begin{equation}
\g^r = \bo_{n_{\A}(\alpha) = r} \g_{\alpha}
\label{gr = sum n(alpha)=r}
\end{equation}
In particular, for $\beta \in \A$, $n_{\A}(\beta)=1$ so $\g_{\beta} \subset \g^1 \subset \m^+$.

Conversely, let $\alpha \in \mathcal R$ such that $n_{\A}(\alpha)=1$. Then there exists a simple root $\beta \in \A$ such that $\alpha = \alpha_0 + \beta$ where $n_{\A}(\alpha_0) = 0$. Consequently,
\[ \g_{\alpha}=[\g_{\alpha_0},\g_{\beta}] \]
where $\g_{\alpha_0} \subset \k_{\A}^{\CM} \subset \h^{\CM}$ and $\g_{\beta} \subset \m^+$. But $\m^+$ represents an invariant almost complex structure on $G/H$ thus it is ${\rm Ad}(H)$-invariant. Then $\g_{\alpha} \subset [\h^{\CM},\m^+] \subset \m^+$ for all $\alpha$ such that $n_{\A}(\alpha)=1$.
\end{proof}

\section{Examples}

In this section we shall number the positive roots $\mathcal B = \{ \beta_1,\beta_2,\ldots \}$ and replace $\A$ in $\p_{\A}$, $\k_{\A}$ or $\l^+_{\A}$ by the corresponding number : for example $\p_{\{ \beta_1,\beta_3,\beta_4} \}$ becomes $\p_{134}$. Also, for a root $\alpha$ written $(a_1,\ldots,a_k)$, in standard coordinates, we abreviate $\g_{\alpha}$ into $[a_1,\ldots,a_k]$ and $\leftr \g_{\alpha} \rightr$ into $\leftr a_1,\ldots,a_k \rightr$.

\subsection{Twistor space of a Wolf space}

Here, we treat a particular case where $K=H$. Let $\t \subset \g$ be a Cartan subalgebra and let $\mathcal R$ be the corresponding root system. We choose a scalar product $\la \ , \ \ra$ on $\t^*$, invariant by the Weyl group (the group of automorphisms of $T$ coming from inner automorphisms of $G$). We define the following quantity, for $\alpha, \beta \in \mathcal R$ : 
\[ (\alpha / \beta) = 2\frac{\la \alpha,\beta \ra}{\la \beta,\beta \ra} \]
According to \cite{bro}, p198, this is an integer, ranging between -3 and 3, that gives information on the relative positions of $\alpha$, $\beta$. Now, let $\gamma$ be the highest weight, for a choice of simple roots $\mathcal B$ : $\forall \alpha \in \mathcal R$, $\la \gamma, \gamma \ra \geq \la \alpha,\gamma \ra$. Moreover $\la \gamma, \gamma \ra = \la \alpha,\gamma \ra$ if and only if $\alpha = \gamma$. Finally $(\alpha / \gamma)$ takes the values $-1$, $0$ and $1$ if $\alpha \neq \gamma$, and $2$ (resp. $-2$) if $\alpha = \gamma$ (resp. $-\gamma$).

The Cartan subalgebra $\t$ is identified with its dual via the scalar product $\la \ , \ \ra$ so that we can see $\gamma$ as an element of $\t$. Let $K$ be the centralizer in $G$ of the 1-dimensional abelian subalgebra generated by $\gamma$. With the description made of the parabolic subalgebras and centralizers of torus (cf theorem \ref{pA} and lemma \ref{k=kA}), this corresponds to
\[ \A = \{ \alpha \in \mathcal B \ | \ (\alpha/\gamma) = 1 \}
\]
and so
\[ \k = \k_{\A} = \t \oplus \bigoplus_{(\alpha / \gamma) = 0} \leftr \g_{\alpha} \rightr \]
We calculate the first canonical series associated to $\p_{\A}$ :
\[ \l^1 = \l_{\A}^+ = \bo_{(\alpha / \gamma) > 0} \g_{\alpha} = \g_{\gamma} \oplus \bo_{(\alpha / \gamma) = 1} \g_{\alpha}, \qquad \l^2 = \g_{\gamma}, \qquad \l^3 = \{0\}. \]
Indeed, for $\alpha,\alpha'$ such that $(\alpha / \gamma) = (\alpha' / \gamma) = 1$, $(\alpha + \alpha'/ \gamma) = 2$ so $\alpha + \alpha' \notin \mathcal R$ or $\alpha + \alpha' = \gamma$. Then 
\[ \g_0 = \k^{\CM}, \quad \g_1 = \bo_{(\alpha / \gamma) = 1} \g_{\alpha}, \quad \g_2 = \g_{\gamma}, \quad \g_3 = \{ 0 \} \]
and the canonical fibration of Burstall and Rawnsley (the "fibration of degree 2") is
\[ G^{\CM}/P_{\A} \to G/F \]
where $F$ is the subgroup of $G$ with Lie algebra
\[ \mathfrak f = \bo_{i \geq 0} \leftr \g_{2i} \rightr = \k \oplus \leftr \g_{\gamma} \rightr \]
Thus, the fibre is 2-dimensional, isomorphic to $\CM P(1)$ (see p132 of \cite{sa} for details), inducing a quaternionic structure on $G/F$ and $G/K$ is the twistor space of a Wolf space that is, of a quaternion-Kähler, compact, symmetric space.

As for the fibration of degree 3, $\g_3 = \{0\}$ so
\[ \h = \leftr \g_0 \rightr = \k. \]
Consequently, $G/K = G/H$ is equipped not only with a Kählerian, but also with a 3-symmetric strictly nearly Kähler structure given by
\[ \m^+ = \g_{-\gamma} \oplus \bo_{(\alpha / \gamma) = 1} \g_{\alpha} \]
This almost complex structure $J_2$ is obtained from the complex structure $J_1$, represented by $\l^+$, by reversing the sign of $J_1$ along the fibre of $G/K \to G/F$. Indeed, the vertical subspace is identified with $\leftr \g_{\gamma} \rightr$. This is another proof, for a Wolf space, of theorem \ref{QK twistor is NK} but also of proposition \ref{sym twistor is 3sym}.

\subsection{$G_2$-spaces}

Take $G=G_2$. The roots may be written
\[
\begin{array}{ll}
(0,2\sqrt 3) & (0,-2\sqrt 3) \\
(3,\sqrt 3)  & (-3,-\sqrt 3) \\
(-3,\sqrt 3) & (3,-\sqrt 3) \\
(1,\sqrt 3)  & (-1,-\sqrt 3) \\
(-1,\sqrt 3) & (1,-\sqrt 3) \\
(2,0)        & (-2,0)
\end{array}
\]
Make the choice of simple roots $\beta_1 = (-3,\sqrt 3)$, $\beta_2=(2,0)$. Then, the roots in the left row are the positive roots :
\begin{eqnarray*}
(-1,\sqrt 3) & = & (-3,\sqrt 3) + (2,0) \\
(1,\sqrt 3)  & = & (-3,\sqrt 3) + 2(2,0) \\
(3,\sqrt 3)  & = & (-3,\sqrt 3) + 3(2,0)
\end{eqnarray*}
and
\[ \gamma = (0,2\sqrt 3) = 2(-3,\sqrt 3) + 3(2,0) \]
is the highest weight. We already know two $G_2$-3-symmetric spaces (and by the classification of Gray and Wolf, see \cite{wo} or \cite{gr2}, these are the only ones). Firstly, by the previous example, the twistor space $G_2/U(2)$ of the exceptionnal Wolf space $G_2/SO(4)$ is Kähler but also nearly Kähler, 3-symmetric. Secondly, the 6-dimensional sphere $S^6 \simeq G_2/SU(3)$ is also a 3-symmetric space. There are 3 types of complex flag manifolds of the form $G_2^{\CM}/P$, corresponding to
\begin{eqnarray*}
\k_1 & = & 2 \RM \oplus \leftr 2,0 \rightr \\
\k_2 & = & 2 \RM \oplus \leftr -3,\sqrt 3 \rightr \\
\k_{12} & = & 2 \RM
\end{eqnarray*}
The first two, $G_2^{\CM}/P_1$ and $G_2^{\CM}/P_2$ constitute two distinct realizations of $G_2/U(2)$ as a {\it complex} flag manifold. The first one corresponds to the twistor space of the exceptionnal Wolf space because $\la (2,0),(0,2\sqrt 3) \ra = 0$, for the canonical scalar product, so $K_1$ is the centralizer of the torus generated by $\gamma$. For the second one, we have
\[ \g^1 = [-1,\sqrt 3] \oplus [2,0], \quad \g^2 = [1,\sqrt 3], \quad \g^3 = [0,2\sqrt 3] \oplus [3,\sqrt 3], \quad \g^4 = \{0\} \]
(we can first compute $\l^i$, $i=1,2,3,4$ and then use the definition \ref{canonical series} of $\g^i$ or apply (\ref{gr = sum n(alpha)=r})). Thus,
\[ \h = 2\RM \oplus \leftr -3,\sqrt 3 \rightr \oplus \leftr 0,2\sqrt 3 \rightr \oplus \leftr 3,\sqrt 3 \rightr, \quad \m^+ = [-1,\sqrt 3] \oplus [2,0] \oplus [-1,-\sqrt 3] \]

Finally, $G/H = G_2/SU(3) \simeq S^6$. Indeed, we showed in \cite{bu2} that the 6-dimensional sphere admits an integrable twistor space with fibre isomorphic to $\CM P(2)$. This is the "reduced twistor space" ${\boldsymbol{\mathcal Y}}$ described above or more precisely one of its connected components. The two components ${\boldsymbol{\mathcal Y}}_1$ and ${\boldsymbol{\mathcal Y}}_2$, see (\ref{Grassmannian bundles}), are isomorphic and the fibrations $G_2/U(2) \simeq {\boldsymbol{\mathcal Y}}_1 \to S^6$ and $G_2/U(2) \simeq {\boldsymbol{\mathcal Y}}_2 \to S^6$ are obtained one from the other by composing with the antipodal map of $\ZZ$, $j \mapsto -j$ or in the present setting by exchanging the roles of $\m^+$ and $\m^-$ in the definition of the fibration of degree 3. Moreover, we can not construct another twistor space on $S^6$ with our method because $SU(3)$ acts transitively on the fibres of ${\boldsymbol{\mathcal Y}}_1$, ${\boldsymbol{\mathcal Y}}_2$.

What is the fibration of degree 3 associated with $\p_{12}$, with total space $G_2/S^1 \times S^1$ ? This time, the second canonical fibration is
\[ \g^1 = [2,0] \oplus [-3,\sqrt 3], \quad \g^2 = [-1,\sqrt 3], \quad \g^3 = [1,\sqrt 3] \]
\[  \g^4 = [3,\sqrt 3], \quad \g^5 = [0,2\sqrt 3], \quad \g^6=\{0\} \]
So
\[ \h = 2\RM \oplus \leftr 1,\sqrt 3 \rightr, \quad \m^+ = [2,0] \oplus [-3,\sqrt 3] \oplus [1,-\sqrt 3] \oplus [3,\sqrt 3] \oplus [0,-2\sqrt 3] \]

The resulting 3-symmetric space is $G_2^{\CM}/P_1 \simeq G_2/U(2)$, the twistor space of $G_2/SO(4)$. To see this, we must make another choice of simple roots : $\beta'_1=(0,-2\sqrt 3), \beta'_2=(1,\sqrt 3)$. Then
\[ (3,-\sqrt 3) = 2(0,-2\sqrt 3) + 3(1,\sqrt 3) \]
is the highest weight, $(1,\sqrt 3)$ satisfies $\la (1,\sqrt 3),(3,-\sqrt 3) \ra=0$ and all the roots whose weight spaces compose $\m^+$, except $(-3,\sqrt 3)=-(3,-\sqrt 3)$, are positive roots as required.

\subsection{Special features associated with $SU(4)$}

The roots of $SU(4)$ are
\[
\begin{array}{ll}
(1,0,0,-1) & (-1,0,0,1) \\
(1,0,-1,0) & (-1,0,1,0) \\
(1,-1,0,0) & (-1,1,0,0) \\
(0,1,0,-1) & (0,-1,0,1) \\
(0,1,-1,0) & (0,-1,1,0) \\
(0,0,1,-1) & (0,0,-1,1)
\end{array}
\]
Let $\beta_1 = (1,-1,0,0)$, $\beta_2 = (0,1,-1,0)$, $\beta_3=(0,0,1,-1)$ be the simple roots. The other positive roots are $(1,0,-1,0)=\beta_1 + \beta_2$, $(0,1,0,-1)=\beta_2 + \beta_3$ and the highest weight $\gamma=(1,0,0,-1)=\beta_1 + \beta_2 + \beta_3$. 

There are $7 = 2^3 -1$ types of flag manifolds built up from $SU(4)$. The associated fibrations of degree 3 are all trivial, in the sense that, $\forall \A \subset \mathcal B$, $H=K_{\A}$ (we will not reproduce all the calculations because they are quite lengthy), except when $\A = \mathcal B$ where the fibration of degree 3 is
\begin{multline}
SU(4)^{\CM}/P_{123} \simeq SU(4)/S^1 \times S^1 \times S^1 \\ \longrightarrow SU(4)^{\CM}/P_{13} \simeq SU(4)/S(U(1) \times U(1) \times U(2))
\label{3-fibration SU(4)/T}
\end{multline}
The last space is the twistor space of the quaternion-Kähler Grassmannian of complex planes in $\CM^4$,  (as in the previous example with $G_2/S^1 \times S^1$, it will appear when changing the basis of the root system into $\beta'_1=(-1,0,1,0)$, $\beta'_2=(1,0,0,-1)$, $\beta'_3=(0,-1,0,1)$). 

In the other cases, we obtain one, and perhaps several different 3-symmetric structures on the flag manifold $SU(4)/K_{\A}$. For example, $SU(4)/S(U(1) \times U(1) \times U(2))$ has three different realizations as a complex flag manifold : $SU(4)^{\CM}/P_{13}$, but also $SU(4)^{\CM}/P_{12}$ and $SU(4)^{\CM}/P_{23}$. As we said, the first one is the twistor space of a Wolf space, $Gr_2(\CM^4)$ (indeed $K_{13}$ is the centralizer of $\gamma$). The others are more interesting because the canonical fibration has rank 4, the fibre being isomorphic to $\CM P(2)$, so the base is a symmetric space but not a Wolf space. Indeed this is $\CM P(3)$. Nevertheless the fibration still has complex, totally geodesic fibres thus the conditions are satisfied to perform the same modification of the Kähler structure as in the quaternion-Kähler situation by \cite{na}. Consider, for example $SU(4)^{\CM}/P_{12}$ (the other case is similar). We calculate the second canonical series :
\[ \g^1 = [1,-1,0,0] \oplus [0,1,-1,0] \oplus [0,1,0,-1], \quad \g^2 = [1,0,0,-1] \oplus [1,0,-1,0], \quad \g^3 = \{0\} \]
Consequently the canonical fibration is $SU(4)/K_{12} \to SU(4)/F$ where the Lie algebra of $F$ is
\[
\mathfrak f = \k_{12} \oplus \leftr 1,0,0,-1 \rightr \oplus \leftr 1,0,-1,0 \rightr 
\]
Thus, $F = S(U(1) \times U(3))$ and the symmetric space $SU(4)/F$ is isomorphic to $\CM P(3)$. Now the fibration {\it of degree 3}, $SU(4)/K_{12} \to SU(4)/H$ is characterized by $H=K_{12}$ and
\[ \m^+ =  [1,-1,0,0] \oplus [0,1,-1,0] \oplus [0,1,0,-1] \oplus [-1,0,0,1] \oplus [-1,0,1,0] \]
so $\m^+$ is obtained from $\l^+_{12}$ by changing the sign along the vertical subspace of the canonical fibration, represented by $\leftr 1,0,0,-1 \rightr \oplus \leftr 1,0,-1,0 \rightr$.

However this 3-symmetric structure coincides with the 3-symmetric structure on the twistor space of a quaternion-Kähler manifold. To see this we must take the basis $\beta''_1=(-1,0,1,0)$, $\beta''_2=(0,0,-1,1)$ and $\beta''_3=(0,1,0,-1)$, so that the highest weight is $\gamma''=(-1,1,0,0)$. The centralizer of the torus generated by $\gamma''$ has Lie algebra
\[ \k''_{13} = 3\RM \oplus \leftr 0,0,1,-1 \rightr = \k_{12}, \]
\[ \l''^+_{13} = [-1,1,0,0] \oplus [0,1,-1,0] \oplus [0,1,0,-1] \oplus [-1,0,0,1] \oplus [-1,0,1,0] \]
and finally the 3-symmetric structure, as in example 1, is represented by
\[ \m''^+ = [1,-1,0,0] \oplus [0,1,-1,0] \oplus [0,1,0,-1] \oplus [-1,0,0,1] \oplus [-1,0,1,0] = \m^+ \]
So even if $M = SU(4)/S(U(1) \times U(1) \times U(2))$ has three distinct invariant complex structures, corresponding to three distinct realizations as a complex flag manifold, or equivalently to three twistor fibrations on a symmetric space with total space $M$, the corresponding 3-symmetric nearly Kähler structures coïncide. In other words, there are several $J_1$ but only one $J_2$, in terms of proposition \ref{QK twistor is NK}. This fact was already noticed by Salamon in \cite{sa2}, section 6.

The other feature associated with $SU(4)$ we found interesting is a kind of counter example of proposition \ref{reciproque}, i.e. a complex twistor space on a 3-symmetric space, not coming from a fibration of degree 3. Again, the 3-symmetric space we consider is the twistor space of a Wolf space, here $Gr_2(\CM^4)$, with the corresponding 3-symmetric structure as in paragraph 6.1 given by :
\[ \h = 3\RM \oplus \leftr 0,1,-1,0 \rightr \]
\[ \m^+ = [-1,0,0,1] \oplus [1,0,-1,0] \oplus [1,-1,0,0] \oplus [0,1,0,-1] \oplus [0,0,1,-1] \]
In view of proposition \ref{twistor space associated to B}, we make a change of base $\mathcal B \mapsto \mathcal B'' = \{ \beta''_1, \beta''_2, \beta''_3 \}$ as above. Then the positive roots are :
\[ (-1,0,1,0) \quad (0,0,-1,1) \quad (0,1,0,-1) \quad (-1,0,0,1) \quad (0,1,-1,0) \quad (-1,1,0,0)    \]
We compute $\n^+$ as in (\ref{n+ associated to B}) :
\[ \n^+ = [-1,0,0,1] \oplus [-1,0,1,0] \oplus [-1,1,0,0] \oplus [0,1,0,-1] \oplus [0,0,-1,1] \]
We remark that $\n^+$ contains all the weight spaces associated to a simple root $\beta''_i$. Thus
\[ [\n^+,\n^+] + \n^+ = \l''^+_{123} \]
In particular $\n^+$ satisfies (\ref{Jt integrable homogene}) so the $H$-twistor space associated to $\n^+$ is integrable, isomorphic to $SU(4)/(S^1)^3 \simeq SU(4)^{\CM}/P_{123}$. However, the first term of the second canonical series associated to $\l''^+$,
\[ \g''^1 = [-1,0,1,0] \oplus \oplus [0,1,0,-1] \oplus [0,0,-1,1], \]
is not included in $\m^+$ (nor in $\m^-$) thus it can't be the fibration of degree 3 by our remark \ref{g1,m+}. Yet, the total space and the base $M \simeq SU(4)/S(U(1) \times U(1) \times U(2))$, equipped with its nearly Kähler structure, are the same. Only the inclusion of $SU(4)/(S^1)^3$ into $\ZZ$, the Riemannian twistor space of $M$, changes.

\subsection{The complex projective space $\CM P(2q+1)$}

The odd dimensional projective space appears, in this subsection, as
\begin{equation}
\CM P(2q+1) \simeq Sp(q+1)/S^1 \times Sp(q) 
\label{CP(2q+1)}
\end{equation}

The properties of this homogeneous realization are very different from those attached to the other realization :
\begin{equation} 
\CM P(n) \simeq SU(n+1)/S(U(1) \times U(n)),
\label{CP(n)}
\end{equation}
appearing in the previous example. For example, it is from (\ref{CP(2q+1)}) that $\CM P(2q+1)$ reveals the twistor space of $\HM P(q) \simeq Sp(q+1)/Sp(1)Sp(q)$. On the other hand, it can be seen on (\ref{CP(n)}) that $\CM P(n)$ is a {\it symmetric} space.

Among all the Wolf spaces, the interest of $\HM P(q)$ is that its Riemannian holonomy is {\it generic}, equal to $Sp(1)Sp(q)$ so it can be used as a model for the study of non locally symmetric, irreducible, quaternion-Kähler  manifolds. Then, the holonomy of $\nb$ on $\CM P(2q+1)$, equal to $S^1 \times Sp(q)$ (recall \ref{connexion normale = intrinseque} that the normal connection coincides with the canonical Hermitian connection on a 3-symmetric space) is also generic among the nearly Kähler manifolds constructed on the twistor space of a quaternion-Kähler manifold with positive scalar curvature.

To simplify we will set $q=2$, the situation for $q=1$ being perhaps too singular. The other cases $q \geq 3$ easily deduce from this one. 

The roots of $Sp(3)$ are, in standard coordinates,
\[
\begin{array}{llllll}
(2,0,0) & (0,2,0) & (0,0,2) & (1,-1,0) & (1,0,-1) & (1,0,-1) \\
(1,1,0) & (1,0,1) & (0,1,1) & (-2,0,0) & (0,-2,0) & (0,0,-2) \\
(-1,1,0)& (-1,0,1)& (-1,0,1)&(-1,-1,0) & (-1,0,-1)& (0,-1,-1)
\end{array}
\]
We choose $\beta_1=(1,-1,0)$, $\beta_2=(0,1,-1)$ and $\beta_3=(0,0,2)$ so that $\gamma=(2,0,0)$. The twistor space of $\HM P(3)$ is $\CM P(5) \simeq Sp(3)^{\CM}/P_{23}$ with complex structure
\[ \l^+ = [2,0,0] \oplus [1,-1,0] \oplus [1,0,-1] \oplus [1,1,0] \oplus [1,0,1], \]
and 3-symmetric structure
\[ \m^+ = [-2,0,0]\oplus [1,-1,0] \oplus [1,0,-1] \oplus [1,1,0] \oplus [1,0,1]. \]
In general we would have, for $\CM P(2q'+1)$,
\begin{multline}
\l^+ = [2,0,\ldots,0] \oplus [1,-1,0,\ldots,0] \oplus [1,0,-1,0,\ldots,0] \oplus \ldots \oplus [1,0,\ldots,0,-1] \\ \oplus [1,1,0,\ldots,0] \oplus [1,0,1,0,\ldots,0] \oplus \ldots \oplus [1,0,\ldots,0,1]
\label{l+ CP(2q+1)}
\end{multline}

Let
\[ \n^+ = [2,0,0] \oplus [-1,1,0] \oplus [1,0,-1] \oplus [1,1,0] \oplus [1,0,1] \]
We have
\[ [\n^+,\n^+]^{\m} = [2,0,0] \oplus [1,1,0], \quad [\n^+,\n^+]^{\h} = [0,2,0] \oplus [0,1,-1] \oplus [0,1,1] \]
so $\n^+$ satisfies
\[ [\n^+,\n^+]^{\m} \subset \n^+, \quad \text{and} \quad [\,[\n^+,\n^+]^{\h},\n^+] = [1,1,0] \subset \n^+ \]

It is not difficult to remark that $[\n^+,\n^+]^{\h}$ is isomorphic to (\ref{l+ CP(2q+1)}) for $q'=q-1=1$ (it suffices to remove the zero in the first place). Since $\h = \k \oplus [\n^+,\n^+]^{\h}$ we have $H/K \simeq \CM P(3)$. This observation is to be made for all $q$. In conclusion, $\forall q \geq 1$, $\CM P(2q+1)$ has a complex twistor space with fibre isomorphic to $\CM P(2q-1)$.

\section{Twistors spaces of 3-symmetric spaces}

The classification of Gray and Wolf \cite{wo} differentiates between three types of compact, simply connected, 3-symmetric spaces $M \simeq G/H$.

\ni (\romannumeral 1) First, the case where $H$ is already the centralizer of a torus. Then $\mathrm{rank}\, G = \mathrm{rank} \, H$ and we can write $\h = \k_{\A}$. More precisely, there are two subcases : \\
\ni a) $\A$ is a singleton : $\A = \{ \delta \}$ where $n_{\delta}(\gamma)=1$ or $2$. \\
\ni b) $\A$ has two distinct elements $\delta_1, \delta_2$ with $n_{\delta_1}(\gamma)=n_{\delta_1}(\gamma)=1$ \\
Then the construction of Burstall, Rawnsley \cite{bur} applies and $M$ is the twistor space of a symmetric space.

\ni (\romannumeral 2) The two groups $G$ and $H$ have same rank but $H$ is not the centralizer of a torus. Then, by \cite{wo}, $H$ is a maximal subgroup of $G$ and there exists a simple root $\delta$ such that $n_{\delta}(\gamma)=3$ and $(\mathcal B - \{\delta\}) \cup \{ -\gamma \}$ constitutes a simple root set for $\h$ which thus equals 
\begin{equation}
\h_{\delta} = \t \oplus \bo_{\substack{n_{\delta}(\alpha) = 0 \\ \ \alpha \in \mathcal R^+}} \leftr \g_{\alpha} \rightr \oplus \bo_{n_{\delta}(\alpha) = 3}  \leftr \g_{\alpha} \rightr 
\label{h isotropy irr.}
\end{equation}

\ni (\romannumeral 3) The third case corresponds to $\mathrm{rank}\, H < \mathrm{rank}\, G$. This includes two exceptionnal homogeneous spaces $Spin(8)/G_2$ and $Spin(8)/SU(3)$ and an infinite family defined by $G = H \times H \times H$ and $H$ stands in $G$ as the diagonal subgroup, the fixed point set of the permutation of order 3.

Thus our assumption that the automorphism of order 3 is inner, or equivalently $\mathrm{rank} \ G = \mathrm{rank} \ H$, appears not too restrictive, since only the spaces in the last group fail to satisfy it.

It is not hard to understand with our formalism what is the 3-symmetric structure on the flag manifolds in (\romannumeral 1). Recall (\ref{gr = sum n(alpha)=r}) : ${\displaystyle \g^r = \bo_{n_{\A}(\alpha) = r} \g_{\alpha}}$.
But $n_{\A}(\gamma) = 1$ or $2$, for the highest weight $\gamma$, so we can't have $n_{\A}(\alpha) \geq 3$. Consequently $\g^3 = \{0\}$, $H=K$ i.e. the fibration of degree 3 is trivial and $G/H$ is itself equipped with a 3-symmetric structure that consists in changing the sign of $\l^+$ along the vertical subspace, represented by $\leftr g^2 \rightr$ (see examples 1 and 3).

The three classes (\romannumeral 1)-(\romannumeral 3) can be characterized in a geometrical way by the holonomy representation of the intrinsic connection $\nb$. This clue observation has a history that starts with Reyes Carri\'on \cite{re}, carries on with Belgun, Moroianu \cite{be} and ends up with the two articles of Nagy \cite{na,na2}. In our case, since all manifolds are homogeneous and the intrinsic connection of a nearly Kähler 3-symmetric space coïncides with the normal connection by proposition \ref{connexion normale = intrinseque}, this is also the isotropy representation. 

\ni (\romannumeral 1) The decomposition
\begin{equation}
\m = \leftr \g^1 \rightr \oplus \leftr \g^2 \rightr 
\label{vert+horiz 3sym}
\end{equation}
is preserved by $\h$ so $\nb$ preserves the horizontal and vertical distributions. Moreover $s$, the automorphism of order 3, is an inner automorphism, for 3-symmetric spaces of type (\romannumeral 1) or (\romannumeral 2). Consequently, the invariant subspaces of ${\rm Ad}(H)$ are stable by $s_*$ and from (\ref{canonical J}), by the canonical almost complex structure $J$. Finally, the isotropy representation of $G/H$ in this case is {\it complex reducible}  (we see tangent spaces as complex vector spaces, identifying $i$ with $J$ at each point).

\ni (\romannumeral 2) On the other hand, the manifolds in (\romannumeral 2) are exceptional homogeneous spaces known to be isotropy {\it irreducible}. In fact these are the only non symmetric isotropy irreducible homogeneous spaces $G/H$ such that ${\rm rank} \, G = {\rm rank} \, H$ (see corollary 8.13.5 of \cite{wo3}).

\ni (\romannumeral 3) Finally  the isotropy, or the holonomy representation of $\nb$ on the spaces of type (\romannumeral 3) is {\it real reducible} i.e. there exists an invariant distribution $\V$ such that $TM = \V \oplus J\V$.

We want to construct twistor spaces over the manifolds of the first two types such that the fibration $G/K \to G/H$ is a fibration of degree 3. According to section 5, this consists, for (\romannumeral 1), in a change of base such that the new base satisfies \ref{reciproque bis}. Examples are given in section 6. They are $G_2/S^1 \times S^1 \to G_2/U(2)$ and (\ref{3-fibration SU(4)/T}).

As for the manifolds of type (\romannumeral 2), $H$ still has maximal rank so \ref{Z is flag} applies and we must look for $H$-twistor spaces over $M$ among the flag manifolds. In view of the description made of $\h$, there is a natural candidate, $G/K_{\delta}$.
\begin{prop}
Let $M=G/H$ be an isotropy irreducible 3-symmetric space where $G$ is a simple Lie group. There exists a simple root $\delta$ such that $n_{\delta}(\gamma) = 3$ and the Lie algebra of $H$ is (\ref{h isotropy irr.}), i.e. $H = H_{\delta}$. Then $G/H$ is the base of the fibration of degree 3 associated to $\p_{\delta}$. Consequently $M$ has a complex twistor space with fibre $H_{\delta}/K_{\delta}$.
\end{prop}
\begin{proof}
By (\ref{gr = sum n(alpha)=r}),
\[ \h_{\delta} = \k_{\delta} \oplus \leftr \g^3 \rightr, \]
for the second canonical series associated to $\p_{\delta}$. The next term is $\g^4 = \{0\}$ because $n_{\delta}(\gamma)=3$. Finally $\h = \bo_{i \in \NM} \leftr \g^{3i} \rightr$ as required.
\end{proof}

\section{Conclusion}

In the classification mentionned in the introduction (see \cite{na2}), we only treated the case of 3-symmetric spaces. However our study might be of some help for the remaining two cases. Indeed, we can see the odd-dimensional complex projective space $\CM P(2q+1)$ considered in paragraph 6.4 as a model for nearly Kähler manifolds $M$ built on the twistor space of a quaternion-Kähler manifold $N$, just like one uses $\HM P(q)$, the quaternionic projective space, as a model for the geometry of non locally symmetric quaternion Kähler irreducible spaces themselves. This is because their holonomy is $Sp(q)Sp(1)$, the isotropy of $\HM P(q)$, from Berger's classification \cite{berg}. Thus, the holonomy of $\nb$ on $M$ is $S^1 \times Sp(q)$ and equals the isotropy of $\CM P(2q+1)$. The study of $M$ gives informations on the quaternion-Kähler manifold $N$. For example, if $M$ is 3-symmetric, $N$ is a symmetric space. Another interesting fact (see \cite{sa2}) is that minimal surfaces in $N$ are projections of holomorphic curves in $M$ with respect to the nearly Kähler almost complex structure.

In the same manner, we can think of manifolds of the third class (i.e. non locally 3-symmetric 6-dimensional nearly Kähler manifolds) as modelled on $S^6 \simeq G_2/SU(3)$. Indeed, let $M$ be a 6-dimensional complete nearly Kähler manifold. It was shown by Belgun, Moroianu \cite{be} that if the holonomy representation of the intrinsic connection $\nb$ is complex reducible, i.e. if the holonomy group is a subgroup of $U(1) \times U(2)$, $M$ is isomorphic to $\CM P(3)$ or $\mathbb F^3$, the space of complex flags in $\CM^3$, and by Nagy \cite{na2} that if the holonomy representation is {\it real} reducible, $M$ is isomorphic to $S^3 \times S^3$. In particular, $M$ is 3-symmetric in both cases. Thus the holonomy representation of non locally 3-symmetric, 6-dimensional, nearly Kähler manifolds is irreducible. Equivalently, the holonomy group of $\nb$ is $SU(3)$, like for the sphere $S^6$. However, twistor spaces over 6-dimensional nearly Kähler manifolds have no compatible complex structure, except for $S^6$ (see \cite{bu2}). 

It was noticed several times that dimension 6 is crucial, in the study of nearly Kähler manifolds (for instance, nearly Kähler manifolds in dimension 6 are Einstein, see \cite{gr3}, and their cone has holonomy $G_2$, see \cite{ba}). In the context of the present article we make the remark that it is already representative of Gray and Wolf's classification of 3-symmetric spaces revisited in section 7. Indeed it supplies examples for the three classes : $S^6$ is isotropy irreducible, $\CM P(3)$, $\mathbb F^3$ are the twistor spaces of 4-dimensional symmetric spaces $S^4$ and $\CM P(2)$, and finally $S^3 \times S^3 \simeq SU(2) \times SU(2) \times SU(2) / SU(2)$ is a representative of the third class.

\vs

The first example of a generalization of the theory of twistor spaces on four dimensional manifolds involving $G$-structures is due to Salamon \cite{sa3} for $G=Sp(q)Sp(1)$. It already makes place for $G$-manifolds with torsion. Indeed, to admit a complex twistor space with fibre $\CM P(1)$, an $Sp(q)Sp(1)$-manifold need not be quaternion-Kähler, it is enough to be {\it quaternionic}, i.e. 3 only of the 6 components of the intrinsic torsion (the tensor that measures the failure, for the holonomy, to reduce to $G$, or for the $G$-structure to be 1-flat, in the terminology of \cite{br}) vanish. Thus, quaternionic manifolds are $Sp(q)Sp(1)$-manifolds with torsion. However these 3 components represent the intrinsic torsion of the underlying $GL(q,\HM)Sp(1)$-manifold (see \cite{sa}, chapter 9) so "quaternionic" means also that the $GL(q,\HM)Sp(1)$-structure is 1-flat and the existence of a complex twistor space depends again on the existence of a torsion-free connection, which doesn't seem to be the case here.

\vs \vs

\vs

\small\noindent LATP, Université de Provence\\
39, rue Frédéric Joliot-Curie\\
13453 Marseille cedex 13, France\\
\texttt{jbbut@cmi.univ-mrs.fr}

\end{document}